\documentclass[A4paper]{amsart}

\usepackage{enumerate}
\usepackage{amsmath,xspace,amssymb,mathrsfs}
\usepackage{color}
\usepackage{hyperref}

\input xy
\xyoption{all}
\xyoption{2cell}
\UseAllTwocells
\CompileMatrices

\renewcommand{\phi}{\varphi}

\newcommand{\A}{{\ensuremath{\mathscr A}}\xspace}
\newcommand{\B}{{\ensuremath{\mathscr B}}\xspace}
\newcommand{\C}{{\ensuremath{\mathscr C}}\xspace}
\newcommand{\D}{{\ensuremath{\mathscr D}}\xspace}
 \newcommand{\I}{{\ensuremath{\mathscr I}}\xspace}
\newcommand{\K}{{\ensuremath{\mathscr K}}\xspace}

\newcommand{\V}{{\ensuremath{\mathscr V}}\xspace}
 \newcommand{\Set}{\textnormal{\bf Set}\xspace}
\newcommand{\sgpd}{\textnormal{\bf scat}\xspace}
\newcommand{\cat}{\textnormal{\bf cat}\xspace}

\newcommand{\Cat}{\textnormal{\bf Cat}\xspace}
\newcommand{\Catcc}{\ensuremath{\Cat_{\textnormal{\bf cc}}}\xspace}

\newcommand{\twocat}{\textnormal{\bf 2-Cat}\xspace}
\newcommand{\PhiCoc}{\textsf{$\Phi$-Coc}\xspace}

\newcommand{\Mnd}{{\sf Mnd}\xspace}
\newcommand{\EM}{{\sf EM}\xspace}
\newcommand{\KL}{{\sf KL}\xspace}

\newcommand{\Id}{\textnormal{\sf Id}\xspace}
\newcommand{\Comp}{\textnormal{\sf Comp}\xspace}

\newcommand{\Hom}{\textnormal{Hom}\xspace}
\newcommand{\Lax}{\textnormal{Lax}\xspace}
\newcommand{\Lan}{\textnormal{Lan}\xspace}
\newcommand{\Ps}{\textnormal{Ps}\xspace}

\newcommand{\coop}{\ensuremath{^{\textnormal{coop}}}}
\newcommand{\co}{\ensuremath{^{\textnormal{co}}}}
\newcommand{\op}{\ensuremath{^{\textnormal{op}}}}
\newcommand{\lax}{\ensuremath{_{\textnormal{lax}}}}

\newcommand{\dm}{\ensuremath{_{\textnormal{dm}}}}

\newcommand{\dl}[2]{\textnormal{dl}(#1,#2)}

\newcommand{\ot}{\otimes}

\newtheorem{theorem}{Theorem}[section]
\newtheorem{lemma}[theorem]{Lemma}
\newtheorem{proposition}[theorem]{Proposition}
\newtheorem{corollary}[theorem]{Corollary}
\newtheorem{predefinition}[theorem]{Definition}
\newenvironment{definition}{\begin{predefinition}\rm}{\end{predefinition}}
\newtheorem{preremark}[theorem]{Remark}
\newenvironment{remark}{\begin{preremark}\rm}{\end{preremark}}
\newtheorem{preexample}[theorem]{Example}
\newenvironment{example}{\begin{preexample}\rm}{\end{preexample}}
\newtheorem{preaussage}[theorem]{}

\renewcommand{\proof}{\noindent{\sc Proof:}\xspace}
\def\endproof{~\hfill$\Box$\vskip 10pt}
\def\ilift#1{\overrightarrow{#1}}
\def\plift#1{\overleftarrow{#1}}

\newcommand{\two}{\ensuremath{\mathbf{2}\xspace}}

\begin{document}

\title[Weak aspects of the theory of monads]{Idempotent splittings, colimit 
completion, and weak aspects of the theory of monads}
\author{Gabriella B\"ohm}
\address{Research Institute for Particle and Nuclear Physics, Budapest \\
H-1525 Budapest 114, P.O.B. 49 \\
Hungary}
\email{\tt G.Bohm@rmki.kfki.hu}
\author{Stephen Lack}
\address{Department of Mathematics\\
Macquarie University NSW 2109 \\
Australia}
\email{\tt steve.lack@mq.edu.au}
\author{Ross Street}
\address{Department of Mathematics\\
Macquarie University NSW 2109 \\
Australia}
\email{\tt ross.street@mq.edu.au}

%\date{Feb 2011}
\begin{abstract}
We show that some recent constructions in the literature, named `weak'
generalizations, can be systematically treated by passing from 2-categories to
categories enriched in the Cartesian monoidal category of Cauchy complete
categories. 
\end{abstract}
\maketitle

\label{firstpage}
\maketitle

The word ``weak'' has been used in category theory with various meanings.
An early example was the weakening of universal properties to ask only 
for the existence of a map, not the existence and uniqueness; this gives,
for example the notion of a {\em weak limit} of a diagram: a cone through 
which every other cone factorizes. 

Often, but not always, such weak limits arise in situations which are
either explicitly or implicitly homotopical, and even though one might
not have uniqueness of factorizations, all such factorizations might be
homotopic: here we have in mind examples like the  
homotopy category (of spaces) and  the category of projective modules
over a ring $R$. Then again, the notion of homotopy might arise because
one is working not just in a category but in a 2-category (or bicategory).
This case led to a second use of the word, coming
from the theory of higher categories: one speaks of weak categories,
or weak functors, or weak limits, to mean notions where all equations are
expected to hold only ``up to higher homotopy''. For instance one has 
the theory of {\em weak $n$-categories}.

The weakness considered in this paper is somewhat different, and originated
in the theory of Hopf algebras.
In the 1990s, questions arising in various areas of
mathematics and even physics suggested the need for a generalization of Hopf
algebras. In subfactor theory, for example, the description of reducible
inclusions required a new symmetry structure. Meanwhile in topology,
invariants of 3-manifolds were constructed that could not be derived from Hopf 
algebras. In some low dimensional quantum field theories, non-integral valued
quantum dimensions occurred, implying that the internal symmetry could not be
described by a Hopf algebra. These phenomena led to the introduction of
{\em weak Hopf algebras} \cite{BNSz:WHAI} which successfully dealt with all
these questions. Apart from Hopf algebras themselves, examples of weak Hopf
algebras include Hayashi's 
face algebras \cite{Ha:face}, Yamanouchi's (finite dimensional) generalized
Kac algebras \cite{Ya:g_Kac}, Ocneanu's paragroups \cite{Oc:q_grd}, and in
particular finite groupoid algebras and their linear duals. 

Large parts of the theory of Hopf algebras have now been generalized to 
the weak context --- we could point, for instance, to the study of
Galois extensions by weak Hopf algebras in \cite{CaDeG:w_entw, Brz:coring, 
AlAletal:w_Gal}, involving a non-commutative version of principal bundles
of structure {\em groupoids} (as opposed to the structure groups of 
the non-weak theory), and to the study of weak Hopf algebras 
in braided monoidal categories \cite{AlAletal:w_smash, PaSt:w_bimonad}.

A weak Hopf algebra $H$ (over a commutative ring $k$) involves, like a 
Hopf algebra, a $k$-module equipped with a $k$-algebra structure and 
a $k$-coalgebra structure, subject to compatibility conditions. The difference
is in the compatibility conditions: although in a weak Hopf 
algebra the comultiplication $H\to H\otimes H$ will still preserve the  
multiplication, the condition that it preserve the unit is replaced by 
a weaker one; similarly the multiplication
$H\otimes H\to H$  will still preserve the comultiplication but need not 
preserve the counit. 

The axioms of a weak Hopf algebra ensure that its category of representations
is monoidal. However, the tensor product of representations is not their
$k$-module tensor product (as it is in the case of a Hopf algebra) but a
certain $k$-module retract of it, obtained by splitting a certain  
idempotent.

In the last two paragraphs, we have already seen two key aspects
of weak Hopf algebras: the weakening of unit conditions, and the splitting
of idempotents. In fact these two are very tightly related, especially
if we consider identities in categories rather than units in
algebras/monoids. Let $A$ and $B$ be categories, and define a {\em semifunctor}
from $A$ to $B$ to be a morphism of the underlying directed graphs which
preserves composition, but which is not required to preserve identities.
Then semifunctors from $A$ to $B$ are in natural bijection with functors
from $A$ to the category $QB$ obtained from $B$ by freely splitting the 
idempotents. We recall the construction and properties of $QB$ in 
Section~\ref{sect:Q} below. The process of freely splitting idempotents in a 
category is often called {\em Cauchy completion}, because it is the case 
$\V=\Set$ of a general construction on enriched categories which in the case
where \V is Lawvere's category $[0,\infty]$ gives the Cauchy completion of a 
(generalized) metric space \cite{Lawvere-metric}. 

We then go on to develop a whole ``weak world'', parallel to the classical
world. At the risk of oversimplifying somewhat, we could summarize the approach
by saying that any classical notion implemented in a 2-category $\K$
should be applied not to $\K$ itself; rather one should first take the
2-category $Q_*\K$, obtained from $\K$ by taking the Cauchy completion of
the hom categories, and then apply the notion there. For example, as in
Section~\ref{sect:monoid} below, we can regard a (strict) monoidal category
$B$ as a one object 2-category. Performing local Cauchy completion, we obtain
a monoidal structure on the Cauchy completion $QB$ of $B$, and one can now
consider monoids not in $B$ but rather in $QB$. This will be our notion of
``weak monoid''.  

In fact because of the variety of meanings of the epithet weak, we have 
decided not to use it as our general naming device; instead we use
the prefix ``demi-'', so in the case of the previous paragraph, we define
a {\em demimonoid} in a monoidal category $B$ to be a monoid in $QB$. 

We gradually work through various other structures, weakening them as
we go. This includes monads and their algebras in Section~\ref{sect:monads}, 
and limits in Section~\ref{sect:limits}. The most important instance of a
limit for our purposes is that of an Eilenberg-Moore object in Section
\ref{sect:EM}. 
Our ultimate goal in Section \ref{sec:ftm} is to develop a weak version of the
formal theory of monads \cite{ftm, ftm2}, building on the start made in
\cite{Bohm:wftm}. In particular, we see that for a 2-category \K in which
idempotent 2-cells split, the 2-category $\EM^w(\K)$ of \cite{Bohm:wftm}
is the free completion of \K, as a 2-category in which idempotent 2-cells
split, under (bicategorical) Eilenberg-Moore objects: see
Corollary~\ref{cor:wftm}.  

The classical formal theory of monads has applications in Hopf algebra
theory. Recall that a bialgebra (over a field) can be regarded as an
opmonoidal monad in the monoidal category of vector spaces (considered as a
single object bicategory). Because of this, the Eilenberg-Moore category of
algebras (i.e. the category of modules over the bialgebra) is monoidal with
the monoidal structure lifted from the category of vector spaces (so that the
forgetful functor is strict monoidal). Moreover, any monad in the
Eilenberg-Moore category (i.e. module algebra over the bialgebra) induces a 
wreath in the category of vector spaces. The corresponding wreath product is
called the `smash product algebra'. 

In a similar way, our weak version occurs in constructions related to
weak Hopf algebras. Weak bialgebras (again, over a field) are `weak
bimonads' in the monoidal category of vector spaces. Weak bimonads
were studied in \cite{BLS:weak_bimonad}. They are monads equipped with the
additional structure which ensures that the Eilenberg-Moore category is
monoidal such that the forgetful functor possesses a so-called separable
Frobenius monoidal structure in the sense of \cite{Szlach}. In this case the
monoidal structure of the Eilenberg-Moore category is weakly lifted from the
category of vector spaces in the sense discussed in the current paper. 
Moreover, any monad in the Eilenberg-Moore category (i.e. any module algebra
over a weak bialgebra) induces a weak wreath in the category of vector spaces;
the corresponding weak wreath product in the sense of this paper is the weak
smash product algebra.  

\subsection*{Acknowledgments.}
Most results in this paper were announced in SL's talk \cite{www} at the 
conference CT2009
in Cape Town, South Africa. It is our pleasure to thank the organizers for
their excellent job and the audience for their interest and helpful comments.
All authors are grateful for partial support from the Australian
Research Council, project DP0771252. GB acknowledges financial support also of
the Hungarian Scientific Research Fund OTKA, grant no. K68195.

\section{Local Cauchy completion}

In later sections of this paper we will perform some constructions in Cauchy
completions of categories, i.e. categories obtained by freely splitting
idempotents in some category. We start with collecting some results about
Cauchy completions and by illustrating the relevance of this process for our
subject.

\subsection{The Cauchy completion functor $\cat\to \cat$.}\label{sect:Q}
Write \cat for the  category of categories and functors; later 
we shall want to consider this also as a 2-category \Cat. 
Define a semicategory to be a directed graph with an associative composition
(no identities assumed), and a semifunctor as the obvious notion of
homomorphism of semicategories. Thus a category is precisely a semicategory 
with identities, and a functor is precisely an identity-preserving
semifunctor between categories. Write \sgpd for the category of semicategories
and semifunctors. 

The evident forgetful functor $U:\cat\to\sgpd$ of course has a left 
adjoint $F:\sgpd\to\cat$ which freely adjoins identities to a semicategory.
But it also has a right adjoint $R:\sgpd\to\cat$ which picks out all 
``potential identities'' in the form of idempotents. Explicitly, for a 
semicategory $S$, the objects of $RS$ are the idempotents $\sigma:s\to s$ in
$S$, and a morphism $(s,\sigma)\to(t,\tau)$ in $RS$ is a morphism $\phi:s\to
t$ in $S$ with $\tau\phi=\phi=\phi\sigma$. The identity on $(s,\sigma)$ is
just $\sigma$. 

The adjunction $U\dashv R$ induces a monad on \cat.
We write $Q=RU$ for the endofunctor, and $q:1\to Q$ for the unit. 
For any category $B$, the component $q:B\to QB$ of the unit exhibits
$QB$ as the Cauchy completion of $B$: the category obtained by freely
splitting the idempotents of $B$.
The functor $Q=RU$ is right adjoint to $FU$ and so preserves all limits
in \cat. 

Thus for categories $A$ and $B$, a semifunctor from $A$ to $B$ is the same as
a functor from $A$ to $QB$. That is, the Kleisli category of the monad $Q$ on
$\cat$ can be regarded as the category of generalized functors which are no
longer compatible with identity morphisms. More generally, we shall see that
many weak notions can be obtained by first applying $Q$, then considering the
usual notion.  

\subsection{The Cauchy completion 2-functor $\Cat\to \Cat$.}\label{sec:Q2}
Of course there is also a 2-category \Cat of categories, functors,
and natural transformations, and $Q$ extends to a 2-monad
on \Cat. As a 2-functor, $Q$ is not a right adjoint, and does not preserve
2-categorical limits in general, although it does preserve products and
equalizers, and so all conical limits. Of particular importance will be
the fact that $Q$ preserves finite products.

\begin{example}
Let $\two$ be the category consisting of two objects and a single
non-identity arrow between them. There is a 2-categorical limit
$A^\two$ called the $\two$-power of $A$, defined by the universal
property 
$$\Cat(X,A^\two)\cong \Cat(X,A)^\two$$
where the right hand side is just the category of arrows in $\Cat(X,A)$. 
Explicitly, $A^\two$ is just the category of functors from $\two$ to $A$.
The 2-functor $Q$ does not preserve the power $A^\two$ strictly: the canonical
comparison functor $Q(A^\two)\to Q(A)^\two$ is not invertible, although it is
a surjective equivalence. 

Indeed, an object of $Q(A^\two)$  is an arrow $\alpha:a_1\to a_2$ in $A$,
along with idempotents $\alpha_1:a_1\to a_1$ and $\alpha_2:a_2\to a_2$
satisfying $\alpha_2\alpha=\alpha\alpha_1$. An arrow in $Q(A^\two)$ from
$(\alpha,\alpha_1,\alpha_2)$ to $(\beta,\beta_1,\beta_2)$ consists of arrows
$\phi_i:a_i\to b_i$ for $i\in\{1,2\}$, satisfying $\beta\phi_1=\phi_2\alpha$
and $\beta_i\phi_i=\phi_i=\phi_i\alpha_i$ for $i\in\{1,2\}$.

An object of $Q(A)^\two$  consists of idempotents $\alpha_i:a_i\to a_i$
for $i\in\{1,2\}$, and a morphism $\alpha':a_1\to a_2$ satisfying
$\alpha_2\alpha'=\alpha'=\alpha'\alpha_1$; note the extra condition compared
to $Q(A^\two)$. A morphism in $Q(A)^\two$ from $(\alpha',\alpha_1,\alpha_2)$
to $(\beta',\beta_1,\beta_2)$ consists of morphisms $\phi_i:a_i\to b_i$
for $i\in\{1,2\}$ satisfying $\beta'\phi_1=\phi_2\alpha'$ and 
$\beta_i\phi_i=\phi_i=\phi_i\alpha_i$.

The comparison functor $K:Q(A^\two)\to Q(A)^\two$ is the image of the identity
functor under the composite map
$$
\xymatrix{
\Cat(A^\two,A^\two)\ar[r]^-\cong&
\Cat(A^\two,A)^\two\ar[r]^-{Q}&
\Cat(Q(A^\two),Q(A))^\two\ar[r]^-\cong&
\Cat(Q(A^\two),Q(A)^\two).
}
$$
Its effect is not perhaps what one might expect: it sends an object
$(\alpha,\alpha_1,\alpha_2)$ of $Q(A^\two)$ to
$(\alpha\alpha_1,\alpha_1,\alpha_2)$, while it sends a morphism  
$(\phi_1,\phi_2)$ to $(\phi_1,\phi_2)$. (In brief, set
$\alpha'=\alpha\alpha_1$.)  Clearly $K$ is faithful; to see that it is full we 
must check that for a morphism $(\phi_1,\phi_2)$ in $Q(A)^\two$ we have
$\beta\phi_1=\phi_2\alpha$, but
$\beta\phi_1=
\beta\beta_1\phi_1=
\beta'\phi_1=
\phi_2\alpha'=
\phi_2\alpha\alpha_1=
\phi_2\alpha_2\alpha=
\phi_2\alpha$
as required. Thus $K$ is fully faithful;  
it is also clearly surjective on objects, and so an equivalence of 
categories, but is not injective on objects. (For instance, if 
$\alpha_1:a_1\to a_1$ is any non-identity idempotent, let $a_2=a_1$ and
$\alpha_2=\alpha_1$, and then $\alpha=\alpha_1$ and $\alpha=1$ give two
different objects of $Q(A^\two)$ which get sent by $K$ to the same object of
$Q(A)^\two$.) 
\end{example}

\begin{example}
Similarly, for any category $C$, there is a 2-categorical limit
$A^C$ called the $C$-power of $A$, defined by $\Cat(X,A^C)\cong\Cat(X,A)^C$.
In general, $Q$ does not preserve such powers, even up to equivalence.  

This can be seen as follows. The functor $q:A\to QA$ induces a fully faithful
map $\Cat(C,A)\to \Cat(C,QA)$, and $\Cat(C,QA)$ is Cauchy complete since $QA$
is; thus there is an induced fully faithful inclusion  
$Q(\Cat(C,A))\to \Cat(C,QA)$, and this is the canonical comparison map
$Q(A^C)\to(QA)^C$. 
It is an equivalence if and only if every $f:C\to QA$ is a retract of some 
functor of the form $qg:C\to QA$, where $g:C\to A$. 

Let $A$ be the free-living idempotent, consisting of a single object $*$ and 
a single non-identity arrow $e$ satisfying $e^2=e$, and let $C=QA$: this
has objects $1$ and $e$. We shall show that the identity functor $1:QA\to QA$
is not a retract of a functor of the form $qg$, where $g:QA\to A$.  To give a
functor $QA\to A$ is equivalently to give a split idempotent in $A$. But the
only idempotent which splits in $A$ is the identity. Thus $g$ would have to be
the map constant at the unique object $*$ of $A$. Now any retract of $qg$
would have to be defined using retracts of the object $1$ of $QA$, but it has
no non-trivial retracts, and so $qg$ has no non-trivial retracts. In
particular the identity functor is not a retract. 
\end{example}

\subsection{The local Cauchy completion 2-functor $\twocat\to
    \twocat$.} 

Since $Q:\Cat\to\Cat$ preserves finite products, it induces a 2-functor
$Q_*:\twocat\to\twocat$ sending a 2-category \K to the 2-category 
$Q_*\K$ with the same objects, obtained by applying $Q$ to each hom-category.
 
\section{Monads}\label{sect:monads}

Monads $(A,t)$ in a 2-category $\K$ are the same as monoids $t$ in the strict
monoidal category $\K(A,A)$. A monad in the local Cauchy completion $Q_*\K$ is
thus a monoid in $Q_*\K(A,A)=Q(\K(A,A))$. Let us call this a {\em weak monad} or
{\em demimonad} in \K. These  were considered in \cite{Bohm:wftm} using the
explicit description of Proposition~\ref{prop:demimonoid} below; and also in
\cite{Wis:q-adj} where the name `$\eta$-symmetric regular quasi-monad' is used.

The 2-natural transformation $q:1\to Q$ induces a 2-natural transformation
$q_*:1\to Q_*$, whose component at a 2-category \K is the inclusion 2-functor
$q_*:\K\to Q_*\K$. This sends monads to monads, and so shows how we can 
regard ordinary monads in \K as demimonads. 

\subsection{Monoids}\label{sect:monoid}

Let $(B,\ot,i)$ be a monoidal category. Since the 2-functor 
$Q:\Cat\to\Cat$ preserves finite products, it sends monoidal objects to 
monoidal objects, and so we obtain a monoidal category 
$(QB,\ot',qi)$, which we usually call $Q(B,\ot,i)$ or just $QB$. 
Explicitly, the tensor product $\ot'$ on $QB$ is given on objects
by $(b,\rho)\ot'(b',\rho')=(b\ot b',\rho\ot\rho')$.

\begin{definition}
A {\em weak monoid} or {\em demimonoid} in $(B,\ot,i)$ is a monoid in
$Q(B,\ot,i)$.  
\end{definition}

We now make this more explicit as follows. 

\begin{proposition}\label{prop:demimonoid}
A demimonoid in $B$ on an object $b$ consists of the following structure:
\begin{enumerate}[(i)]
\item an associative multiplication $\mu:b\ot b\to b$;
\item a map $\eta:i\to b$ for which 
\begin{itemize}
\item the composite
$$\xymatrix{
i \cong i\ot i \ar[r]^{\eta\ot\eta} & b\ot b \ar[r]^{\mu} & b }$$
is just $\eta$;
\item the composites
$$\xymatrix{
b \cong b\ot i \ar[r]^{b\ot\eta} & b\ot b \ar[r]^{\mu} & b} 
\qquad \textrm{and}\qquad 
\xymatrix{
b \cong i\ot b \ar[r]^{\eta\ot b} & b\ot b \ar[r]^{\mu} & b 
}$$
are equal (let us call them $\mu_1$);
\item the above map $\mu_1:b\to b$ satisfies $\mu_1\mu=\mu$.
\end{itemize}
\end{enumerate}
\end{proposition}

\proof
Given structure as in the proposition, first note that the composite
$\mu_1\mu_1$ is given by 
$$\xymatrix @C=35pt{
b \ar@/^1pc/[rrr]^{\mu_1} \ar[r]|-{\,\cong\,} \ar[drrr]_{\mu_1} & b\ot i
\ar[r]|-{\,b\ot\eta\,} & b\ot b \ar[r]|-{\,\mu\,} \ar[dr]_{\mu} & b
\ar[d]^{\mu_1} \\  
&&& b }$$
and so $\mu_1$ is idempotent. Thus $(b,\mu_1)$ is an object of $QB$.
The first condition in (ii) says that $\mu_1\eta=\eta$, which in turn says that
$\eta$ is a map in $QB$ from $qi$ to $(b,\mu_1)$. Next we show that
$\mu$ is a map in $QB$ from $(b,\mu_1)\ot(b,\mu_1)$ to $(b,\mu_1)$.
One half of this is the last requirement in (ii), the other half
says that  $\mu(\mu_1\ot\mu_1)=\mu$; or, equivalently 
$\mu(\mu_1\ot b)=\mu=\mu(b\ot\mu_1)$. To see the first of
these, note that the diagram
$$\xymatrix{
b\ot b \ar[r]^-{\cong} \ar[d]_{\mu} & 
i\ot b\ot b \ar[r]^{\eta\ot b\ot b} \ar[d]^{i\ot\mu} &
b\ot b\ot b \ar[r]^{\mu\ot b} \ar[d]^{b\ot\mu} & 
b\ot b \ar[d]^{\mu} \\
b \ar[r]_-{\cong} & i\ot b \ar[r]_{\eta\ot b} & b\ot b \ar[r]_{\mu} & b 
}$$
commutes, and that the bottom path is $\mu_1\mu$ which is indeed 
$\mu$; the second equation is proved similarly.

Thus we have maps $\mu:(b,\mu_1)\ot(b,\mu_1)\to(b,\mu_1)$ 
and $\eta:qi\to(b,\mu_1)$ in $QB$. Associativity is part (i), while the unit
laws are the conditions in the second item in (ii), and so these give a monoid
in $QB$.  

Conversely, any monoid in $QB$ involves an underlying object 
$(b,\mu_1)$ with $\mu_1$ idempotent, equipped with morphisms
$\mu:(b,\mu_1)\ot(b,\mu_1)\to(b,\mu_1)$ and 
$\eta:qi\to(b,\mu_1)$. Associativity of $\mu$ gives (i), the
unit laws give the second condition in (ii), the fact that $\mu$ is a morphism
in $QB$  gives the last condition, and the unit laws and the fact that $\eta$
is a morphism in $QB$ give the first condition in (ii).
\endproof

\begin{example}
As any adjunction in a 2-category induces a monad, any adjunction in the local
Cauchy completion $Q_*\K$ of a 2-category $\K$ induces a demimonad (i.e. a
demimonoid in the hom category $\K(A,A)$) which we describe presently. 

A 1-cell in $Q_*\K$ is a pair $(x,\overline x)$ consisting of a 1-cell $x$,
and an idempotent 2-cell $\overline x:x\to x$, in $\K$. An adjunction in
$Q_*\K$ is given by 1-cells $(x,\overline x):A\to B$ and $(y,\overline y):B\to
A$ together with 2-cells $\phi:(x,\overline x)(y,\overline y) \to (1_{B},1)$
and $\psi:(1_{A},1)\to (y,\overline y)(x,\overline x)$ in $Q_*\K$; obeying the
usual triangle identities. Using the triangle conditions  
$$
\overline x=\phi x.x\psi
\qquad \textrm{and}\qquad
\overline y= y\phi.\psi y,
$$
we can express $\overline x$ and $\overline y$ in terms of $\phi$ and
$\psi$. Thus only the normalization conditions $\phi.\overline x
y=\phi=\phi.x\overline y$ and $\overline y x.\psi=\psi=y\overline x.\psi$ are
left. These are equivalent to commutativity of the diagrams
$$
\xymatrix{
xy\ar[r]^-{x\psi y}\ar[rd]_-\phi&xyxy\ar[d]^-{\phi\phi}\\&1_B
}\quad \textrm{and}\quad
\xymatrix{
1_A\ar[r]^-{\psi\psi}\ar[rd]_-\psi&yxyx\ar[d]^-{y\phi x}\\
&yx
}
$$
Summarizing, an adjunction in $Q_*\K$ is given by 1-cells $x:A\to B$ and
$y:B\to A$ together with 2-cells $\phi:xy\to 1_B$ and $\psi:1_A\to yx$ in
$\K$; rendering commutative these two diagrams. This structure is discussed in
\cite{Wis:q-adj} under the name `regular adjunction context'. The
corresponding demimonad is $(A,yx)$ with the associative multiplication
$\overline y\phi \overline x$ and demiunit $\psi$. 
\end{example}

\subsection{Algebras for monads}\label{sect:algebras}

A monad $(A,t)$ in a 2-category $\K$ induces a monad $\K(B,t)$ on the category
$\K(B,A)$ for any object $B$ of $\K$. We may consider its Eilenberg-Moore
algebras; i.e. the actions of $(A,t)$ on morphisms $B\to A$. We can now define
demiactions of our demimonads as ordinary actions in $Q_*\K$ of the monads in
$Q_*\K$. Even if we start with an actual monad in $\K$, this gives something
new. 

In view of Proposition \ref{prop:demimonoid}, a demimonad in a 2-category $\K$
is given by a 1-cell $t:A\to A$, an associative multiplication $\mu:t^2\to t$
and a 2-cell $\eta:A\to t$ subject to the conditions in  Proposition
\ref{prop:demimonoid}. For the idempotent 2-cell $\mu.t\eta=\mu.\eta t:t\to t$
we write $\mu_1$. 

\begin{proposition}\label{prop:demiaction}
A demiaction of a demimonad $(A,t)$ on a morphism $a:B\to A$ 
consists of a 2-cell $\alpha:ta\to a$ satisfying the associative 
law $\alpha.t\alpha=\alpha.\mu a$ as well as $\alpha.\mu_1 a=\alpha$;
when $t$ is a monad, then $\mu_1=1$ and the second condition is 
automatic. 
\end{proposition}

\proof
An action in $Q_*\K$ consists of a morphism $a:B\to A$ equipped with 
an idempotent $\bar{a}:a\to a$, and an action 
$\alpha:(t,\mu_1)(a,\bar{a})\to(a,\bar{a})$. In order for $\alpha$ to be
a morphism in $Q_*\K$, we need $\bar{a}.\alpha=\alpha=\alpha.\mu_1 a.t\bar{a}$;
or equivalently $\bar{a}.\alpha=\alpha=\alpha.\mu_1 a=\alpha.t\bar{a}$.
The associative law says $\alpha.t\alpha=\alpha.\mu a$ and the unit 
law says $\alpha.\eta a=\bar{a}$. Thus $\bar{a}$ can be recovered from
$\alpha$. We must show that any $\alpha:ta\to a$ satisfying
$\alpha.t\alpha=\alpha.\mu a$ and $\alpha.\mu_1 a=\alpha$ satisfies the
remaining conditions.

First of all $\alpha.\eta a.\alpha=\alpha.t\alpha.\eta ta=\alpha.\mu a.\eta ta
=\alpha.\mu_1 a=\alpha$, and so $\bar{a}\alpha=\alpha$. Furthermore
this gives $\bar{a}\bar{a}=\bar a.\alpha.\eta a=\alpha.\eta
a=\bar{a}$ and $\bar{a}$ is idempotent. Finally
$\alpha.t\bar{a}=\alpha.t\alpha.t\eta a=\alpha.\mu a.t\eta a=\alpha.\mu_1
a=\alpha$.  
\endproof

\begin{remark}\label{rmk:iso}
A morphism $(b,\beta)\to(c,\gamma)$ of $t$-(demi)actions is a 2-cell
$\phi:b\to c$ making the following diagrams commute. 
$$\xymatrix{
tb \ar[r]^{t\phi} \ar[d]_{\beta} & tc \ar[d]^{\gamma} \\
b \ar[r]_{\phi} & c}\qquad\qquad
\xymatrix{
b \ar[r]^{\phi} \ar[dr]^{\phi} \ar[d]_{\beta.\eta b} & c \ar[d]^{\gamma.\eta c} \\ 
b \ar[r]_{\phi} & c. 
}$$
Commutativity of the diagram on the left is the usual condition for
morphisms of $t$-actions; as for the diagram on the right, the exterior 
will commute if the diagram on the left does -- the new condition is that
the two equal paths around the exterior are themselves equal to the diagonal.
Of course if $\beta$ and $\gamma$ are genuine (unital) actions, then
$\beta.\eta b$ and $\gamma.\eta c$ are identities, and this is
automatic.  
It is important to note that the identity morphism on a demiaction
$(b,\beta)$ is given by  $\beta.\eta b:b\to b$; this is the identity only in
the case of a genuine action. It follows that a morphism
$\phi:(b,\beta)\to(c,\gamma)$ 
of demiactions can be invertible without $\phi:b\to c$ being invertible
in \K.
\end{remark}

\subsection{The Eilenberg-Moore category}
A proper monad $(A,t)$ in $\K$ can be regarded as a demimonad; i.e. a monad
$(A,t)$ in $Q_*\K$. For any other object $B$ in $\K$, the induced monad
$Q_*\K(B,t)$ on the category $Q_*\K(B,A)$ is equal to the image of the
monad $\K(B,t)$ on $\K(B,A)$ under the 2-functor $Q:\Cat\to \Cat$ in Section
\ref{sec:Q2}. Hence to give a demiaction of the demimonad $(A,t)$ is
equivalently to give an actual algebra of the latter monad on $Q(\K(B,A))$. 

We may apply this reasoning to the particular 2-category $\K=\Cat$ and its
terminal object $B=1$. Then for any monad $t$ on a category $A$, there is a
coinciding notion of a demiaction of the demimonad $(A,t)$ (i.e. action of
the monad $(A,t)$ in $Q_*\Cat$) and that of an actual algebra of the monad
$Qt=(Qt,Q\mu,Q\eta)$ on $QA$. We call it a {\em $t$-demialgebra}. From
Proposition \ref{prop:demiaction} we obtain the following explicit description.

\begin{proposition}
A $t$-demialgebra for a monad $t$ is the same thing as an object $b\in A$
equipped with a morphism $\beta:tb\to b$ satisfying the associative law
$\beta.t\beta= \beta.\mu b$. 
\endproof
\end{proposition}

We obtain a category $A^{(t)}$ of $t$-demialgebras, by taking as morphisms
the $Qt$-morphisms between the corresponding $Qt$-algebras. This gives 
an isomorphism of categories $A^{(t)}\cong(QA)^{Qt}$. Explicitly, if 
$(b,\beta)$ and $(c,\gamma)$ are $t$-demialgebras, a morphism of demialgebras
from $(b,\beta)$ to $(c,\gamma)$ is a morphism $\phi:b\to c$ satisfying
$\phi.\beta=\gamma.t\phi$ and $\gamma.\eta c.\phi=\phi$. 
As said in Remark \ref{rmk:iso}, the identity morphism on a $t$-demialgebra
$(b,\beta)$ is $\beta.\eta b$.  

\begin{example}\label{ex:identity}
A demialgebra for the identity monad on $B$ is a morphism 
$\beta:b=1_B(b)\to b$ satisfying $\beta\beta=\beta$; that is, an idempotent 
in $B$. In symbols: $A^{(1)}=QA$.
\end{example}

\begin{proposition}
A $t$-demialgebra $(b,\beta)$ is isomorphic to a $t$-algebra if and 
only if the idempotent $\beta.\eta b$ splits. 
\end{proposition}

\proof
Let $(a,\alpha)$ be a $t$-algebra. An isomorphism $(a,\alpha)\cong(b,\beta)$
consists of $t$-demialgebra maps $\sigma:(a,\alpha)\to(b,\beta)$ and
$\pi:(b,\beta)\to(a,\alpha)$ satisfying $\pi\sigma=1$ and
$\sigma\pi=\beta.\eta b$. So certainly if $(b,\beta)$ is isomorphic  
to a $t$-algebra then the idempotent splits. Suppose conversely that the
idempotent splits, say as $\beta.\eta b=\sigma\pi$, with $\pi\sigma=1_a$.
Then $a$ inherits a unique demialgebra structure $\alpha:ta\to a$ such that
$\sigma$ and $\pi$ are both demialgebra morphisms. It remains to check that
$(a,\alpha)$ is in fact an algebra. 

Since $\pi:(b,\beta)\to(a,\alpha)$ is a demialgebra morphism, we have
$\alpha.\eta a.\pi=\pi$; but $\pi\sigma=1$ and so $\alpha.\eta a=1$. 
\endproof

\begin{remark}
We saw before that $Q:\Cat\to\Cat$ does not preserve powers; it also 
does not preserve Eilenberg-Moore objects, since the canonical comparison
$Q(C^t) \to (QC)^{Qt}$
is not invertible; indeed this time it is not even an equivalence in general.
It will be an equivalence if and only if each $t$-demialgebra is a 
retract (in the category of $t$-demialgebras) of a $t$-algebra; in particular,
this will be the case if idempotents split in $C$.  

In more detail, an object of $Q(C^t)$ consists of a $t$-algebra $(A,a)$ and
an idempotent $t$-morphism $e:(A,a)\to (A,a)$. A morphism from $(A,a,e)$ to 
$(B,b,d)$ consists of a $t$-morphism $f:(A,a)\to (B,b)$ satisfying $df=f=fe$.

An object of $Q(C)^{Q(t)}$ consists of an object $A\in C$ with a morphism
$a':tA\to A$ satisfying the associative law $a'.ta'=a'.\mu A$. A  morphism
from $(A,a')$ to $(B,b')$ consists of a morphism $f:A\to B$  satisfying
$fa'=b'.tf$ and $fa'.\eta A=f$.  

The comparison functor $K$ sends $(A,a,e)$ to $(A,ea)$ and a morphism 
$f:(A,a,e)\to(B,b,d)$ to $f:(A,ea)\to(B,db)$. Clearly this is faithful;
while given an arbitrary $f:(A,ea)\to(B,db)$ we have
$f=fea.\eta A=fe$ and $f=fea.\eta A=db.tf.\eta A=db.\eta B.f=df$, and so also
$fa=fea=db.tf=b.td.tf=b.t(df)=b.tf$, which proves that $f$ is a 
morphism $(A,a,e)\to(B,b,d)$ and so that $K$ is also full. 

For any object $(B,b,d)$ of $Q(B^t)$, $K(B,b,d)=(B,db)$ is a retract in
$(QB)^{Qt}$ of the $t$-algebra $(B,b)$ via the epimorphism $d:(B,b) \to
(B,db)$ and its section $d:(B,db)\to (B,b)$. Hence a $t$-demialgebra
$(a,\alpha)$ will be isomorphic to an object in the image of $K$ if and only
if it is a retract of a $t$-algebra. Thus $K$ will be an equivalence if every
$t$-demialgebra is a retract (in the category of demialgebras) of a 
$t$-algebra. 
It will of course be an equivalence whenever idempotents split in $C$.
In general, however, $K$ needs not be surjective on objects, or even 
essentially surjective, and so will not be an equivalence of categories. 

As an example,
consider the category of categories with chosen initial object,
and functors preserving the chosen initial object. This has a subcategory
$B$ consisting of the finite ordinals $n=\{0<1<\ldots<n-1\}$ with $n\ge 2$,
and the category \I with a chosen initial object $0$ and another initial  
object $0'$; we include all maps except that we only allow functors
$n\to\I$ which are constant at 0. There is an evident monad $t$ which 
adjoins a top element (except that when applied to \I it first collapses
$0$ and $0'$ to a single element $0$). Each $n$ has a unique $t$-algebra
structure $\alpha:n+1\to n$ which collapses the top two elements of $n+1$.
The unique map $2\to\I$ makes \I into a demialgebra (but not an algebra).
Any map $\phi:\I\to n$ must satisfy $\phi(0)=0$; but to be a demialgebra
map it would also need to satisfy $\phi(0)=n-1$ which is clearly impossible
for $n\ge 2$. Thus there is {\em no} demialgebra map from \I to a $t$-algebra,
and so certainly \I is not a retract of a $t$-algebra.
\end{remark}

\subsection{Monoids as algebras of the free monoid monad}
We shall now work through in some detail a not entirely trivial 
example. Let $B$ be a monoidal category with countable
coproducts over which the tensor product distributes. We
write, for convenience, as if $B$ were strict. Then 
free monoids in $B$ can be constructed via the usual geometric
series $tb=\sum_n b^n$, where $b^n$ denotes the $n$-fold tensor
power of an object $b$. Then $t$ becomes a monad on $B$, cf. \cite[p. 172,
Theorem 2]{CWM}. A $t$-demialgebra consists of an object $(b,\rho)$ of $QB$
equipped with an action $\beta$ of $Qt$. To give a map $\beta:tb\to b$ is to
give a map $\beta_n:b^n\to b$ for each $n$. The unit law $\beta.\eta b=\rho$
says that $\beta_1=\rho$. 
The fact that $\beta$ is a morphism $(tb,t\rho)\to(b,\rho)$ in $QB$
amounts to the equations $\rho\beta_n=\beta_n$ and $\beta_n\rho^n=\beta_n$.
The associativity constraint can be written as commutativity of
\begin{equation}\label{eq:assoc}
\xymatrix{
b^{\sum_{k=1}^n m_k}\ar[rr]^-{\beta_{\sum_{k=1}^n m_k}}
\ar[d]_-{\otimes_{k=1}^n \beta_{m_k}}&&
b\ar@{=}[d]\\
b^n\ar[rr]^-{\beta_n}&&b
}
\end{equation}
for all natural numbers $n, m_1,\dots,m_n$. Putting $(n=2, m_1=0, m_2=1)$ and
$(n=2, m_1=1, m_2=0)$, \eqref{eq:assoc} reduces to the conditions
\begin{equation}\label{eq:iii}
\beta_2(\beta_0\ot b)=\beta_1=\beta_2(b\ot\beta_0).
\end{equation}
Evaluating \eqref{eq:assoc} at $(n=2, m_1=1, m_2=2)$ and $(n=2, m_1=2,
m_2=1)$ gives the associativity of $\beta_2$. Finally, taking $(n=2,
m_1=p-1, m_2=1)$ for any positive integer $p$, and iterating the resulting
relation, we obtain 
$$
\beta_p=\beta_2(\beta_{p-1}\otimes b)=\dots=
\beta_2(\beta_2 \otimes b)\dots(\beta_2\otimes b^{p-2}).
$$
Together with the associativity of $\beta_2$, this identity implies
commutativity of \eqref{eq:assoc} for any values of $n$ and
$m_1,\dots,m_n$.   
Putting all that together, we see that the entire structure consists of (i) an
associative multiplication $\beta_2:b\ot b\to b$, (ii) an idempotent
$\beta_1:b\to b$ satisfying $\beta_1\beta_2=\beta_2=\beta_2(\beta_1\ot
b)=\beta_2(1\ot\beta_1)$, and (iii) a map $\beta_0:i\to b$ satisfying
$\beta_1\beta_0=\beta_0$ and \eqref{eq:iii}. 
Of course $\beta_1$ is determined by $\beta_2$ and $\beta_0$. A morphism of
$t$-demialgebras from $(b,\beta_2,\beta_0)$ to $(c,\gamma_2,\gamma_0)$
is a morphism $\phi:b\to c$ commuting with the structure maps and satisfying
$\phi\beta_1=\phi$. 

Comparing the above description of the category of $t$-demialgebras and the
category of demimonoids in Section \ref{sect:monoid}, we obtain the
following. 

\begin{theorem}
Suppose that $B$ is a monoidal category with countable coproducts,
and that the monoidal structure distributes over the coproducts. Then 
the category of demimonoids in $B$ is just the category $B^{(t)}$ of 
$t$-demialgebras.
\end{theorem}

This extends the well-known isomorphism between the category of monoids in $B$
and the category $B^t$ of algebras for the free monoid monad $t$. 

\subsection{The 2-category of monads}
Given a monad $t$ on an object $A$ of a  2-category \K, an action
of $t$ on a morphism $a:B\to A$ is a special case of the notion of 
morphism of monad. In fact, for every object $B\in\K$, there is 
an identity monad $1$ on $B$, and to give a morphism $a:B\to A$ and
an action of $t$ on $a$ is equivalently to give a monad morphism
from $(B,1)$ to $(A,t)$. Similarly, one has 2-cells between monad
morphisms, and indeed a whole 2-category $\Mnd(\K)$ of monads in \K.
This was introduced and studied in \cite{ftm}; there is also a variant
$\EM(\K)$ with a different notion of 2-cell which was proposed in
\cite{ftm2}. In the subsequent sections we shall develop weak analogues of
these. 

\section{Lax functors, lax natural transformations and modifications}
\label{sect:Ps}

For any 2-categories $\C$ and $\K$, there is a bicategory of lax functors
$\C\to \K$, lax natural transformations between them, and their
modifications. Its `weak' analog is obtained below by replacing the target
2-category $\K$ by its local Cauchy completion $Q_*\K$.

\subsection{Lax functors}
Let \C and \K be 2-categories. 
The notion of lax functor from \C to \K was introduced in \cite{bicategories}
under the name ``morphism of bicategories'', in the more general context where
\C and \K were bicategories. A lax functor $F$ differs from a 2-functor by the
property that it preserves horizontal composition and identity 1-cells only
up-to natural transformations $\mu:F(-).F(-)\to F(-.-)$ and $\eta:1_{F(-)}\to
F(1_{(-)})$, respectively, which obey coherence conditions called
associativity and unitality.  

\begin{example}
Consider the case $\C=1$, where $1$ is the 2-category with a single 
object $*$ and trivial hom-category $1(*,*)=1$. To give the object-part
of a lax functor $1\to\K$ is to give an object $A\in\K$; to give the
functors between hom-categories is to give a functor $1\to\K(A,A)$;
that is, to give a 1-cell $t:A\to A$ in \K. There is only one component
of $\mu$ to worry about: it is a 2-cell $\mu:t^2\to t$. Similarly the
only component of $\eta$ is a 2-cell $\eta:1\to t$. The associativity
and unit conditions say precisely that $(t,\mu,\eta)$ is a monad.
\end{example}

As was already anticipated in the previous section, since 
a monad in \K is a lax functor $1\to\K$, a demimonad in \K is a lax
functor $1\to Q_*\K$. We therefore define, more generally, a
{\em lax demifunctor} from \C to \K to be a lax functor from \C to
$Q_*\K$. 

These lax demifunctors will be of less importance themselves than 
their morphisms, introduced below. Nonetheless we shall take the trouble
to spell out the structure in more direct terms. First of all, for
each object $C\in\C$, an object $FC\in\K$ is given. For all objects
$C,D\in\C$, a functor $\C(C,D)\to Q_*\K(FC,FD)$ is given; that is,
a functor $\C(C,D)\to Q(\K(FC,FD))$, or equivalently, a semifunctor
$\C(C,D)\to\K(FC,FD)$. As usual, an identity 2-cell $1_f$ in \C will
be sent to an idempotent 2-cell $F1_f:Ff\to Ff$ in \K. For 
1-cells $f:C\to D$ and $g:D\to E$ in \C, we have a 2-cell 
$\mu:Fg.Ff\to F(gf)$ in \K. It is natural in $f$ and $g$, and satisfies
the usual associativity condition, expressed by commutativity of the first
diagram below, as well as a normalization condition which states that the two
triangles in the second diagram commute, for all arrows $e:B\to C$, $f:C\to
D$, $g:D\to E$ in \C:
$$
\xymatrix{
Fg.Ff.Fe \ar[r]^{\mu_{g,f}.1_{Fe}} \ar[d]_{1_{Fg}.\mu_{f,e}} & 
F(gf).Fe \ar[d]^{\mu_{gf,e}} \\
Fg.F(fe) \ar[r]_{\mu_{g,fe}} & F(gfe) }\qquad\qquad
\xymatrix{
Fg.Ff \ar[r]^{\mu} \ar[d]_{F1_g.F1_f} \ar[rd]^\mu& F(gf) \ar[d]^{F1_{gf}} \\
Fg.Ff \ar[r]_{\mu} & F(gf) }$$
Similarly, for each $C\in\C$ 
there is a 2-cell $\eta:1_{FC}\to F1_C$, but the unit conditions now
say that the composites
$$
\xymatrix{
Ff \ar@{=}[r] & Ff.1_{FC} \ar[r]^-{1_{Ff}.\eta} & Ff.F1_C
\ar[r]^-{\mu_{f,1_C}} & 
Ff}
\quad \textrm{and}\quad
\xymatrix{
Ff \ar@{=}[r] & 1_{FD}.Ff \ar[r]^-{\eta.1_{Ff}} & F1_D.Ff
\ar[r]^-{\mu_{1_D,f}} & 
Ff 
}$$
are equal to the idempotent $F1_f$; this time the normalization
condition states that the composite
$$\xymatrix{
1_{FC} \ar[r]^{\eta} & F1_C \ar[r]^{F1_{1_C}} & F1_C }$$
is just $\eta$. 

\subsection{Lax natural transformations}
Once again, any lax functor $F:\C\to\K$ determines a lax demifunctor
$q_*F:\C\to Q_*\K$ with which it is identified. For such a lax demifunctor,
$F1_f=1_{Ff}$ for any 1-cell $f$ in $\C$. Even for such 
lax functors $F,G:\C\to\K$, however, we obtain a new type of morphism,
namely the lax natural transformations $q_*F\to q_*G$. 

What then is a lax natural transformation between lax functors $F,G:\C\to Q_*
\K$? For each $C\in\C$ we should give a 1-cell $FC\to GC$ in $Q_*\K$;
in other words, a 1-cell $xC:FC\to GC$ along with an idempotent
2-cell $\overline{x}C:xC\to xC$. Next, for each 1-cell $f:C\to D$ 
in \C, we should give a 2-cell
$(Gf,G1_f)(xC,\overline{x}C)\to(xD,\overline{x}D)(Ff,F1_f)$ in $Q_*\K$;
that is, a 2-cell 
$$\xymatrix @R1pc {
FC \ar[r]^{xC} \ar[dd]_{Ff} & GC \ar[dd]^{Gf} \\
\relax\rtwocell<\omit>{\ \ xf} & {} \\
FD \ar[r]_{xD} & GD }$$
for which the two composites
$$\xymatrix @R1pc {
Gf.xC \ar[r]^{xf} & xD.Ff \ar[r]^{\overline{x}D.F1_f} & xD.Ff}
\quad \textrm{and}\quad
\xymatrix @R1pc {
Gf.xC \ar[r]^{G1_f.\overline{x}C} & Gf.xC \ar[r]^{xf} & xD.Ff}$$
are both just $xf$. 
This $x$ obeys the same naturality condition and the same compatibility with
composition as a usual lax natural transformation between lax functors $\C\to
\K$: the same diagrams

\begin{minipage}{0.4\linewidth}
  \begin{equation}
    \label{eq:LN0}\tag{LN0}
\xymatrix{
Gf.xC \ar[r]^{xf} \ar[d]_{G\alpha.1_{xC}} & xD.Ff \ar[d]^{1_{xD}.F\alpha} \\
Gf'.xC \ar[r]_{xf'} & xD.Ff' }    
  \end{equation}
\end{minipage}
\begin{minipage}{0.6\linewidth}
  \begin{equation}
    \label{eq:LN1}\tag{LN1}
\xymatrix{
Gg.Gf.xC \ar[r]^{1_{Gg}.xf} \ar[d]_{\mu^G_{g,f}.1_{xC}} & 
Gg.xD.Ff \ar[r]^{xg.1_{Ff}} & xE.Fg.Ff \ar[d]^{1_{xE}.\mu^F_{g,f}} \\
G(gf).xC \ar[rr]_{x(gf)} && xE.F(gf) }
  \end{equation}
\end{minipage}

% $$
% \mathrm{(LN0)}
% \xymatrix{
% Gf.xC \ar[r]^{xf} \ar[d]_{G\alpha.1_{xC}} & xD.Ff \ar[d]^{1_{xD}.F\alpha} \\
% Gf'.xC \ar[r]_{xf'} & xD.Ff' }
% \qquad
% \mathrm{(LN1)}
% \xymatrix{
% Gg.Gf.xC \ar[r]^{1_{Gg}.xf} \ar[d]_{\mu^G_{g,f}.1_{xC}} & 
% Gg.xD.Ff \ar[r]^{xg.1_{Ff}} & xE.Fg.Ff \ar[d]^{1_{xE}.\mu^F_{g,f}} \\
% G(gf).xC \ar[rr]_{x(gf)} && xE.F(gf) }
% $$
\noindent
commute for all 1-cells $f,f':C\to D$, $g:D\to E$ and 2-cells $\alpha:f\to
f'$. The third condition, expressing compatibility with identities, is changed
because the identity 2-cell in $Q_*\K$ on $(xC,\overline{x}C)$ is
$\overline{x}C$. Thus the new condition becomes commutativity of 
\begin{equation}\label{eq:DLN2}
\tag{DLN2}
\xymatrix{
1_{GC}.xC \ar[d]_{\eta^G.\overline{x}C} \ar@{=}[r] & 
xC.1_{FC} \ar[d]^{\overline{x}C.\eta^F} \\
G1_C.xC \ar[r]_{x1_C} & xC.F1_C. }
\end{equation}
For lax functors $F,G:\C\to \K$, we may consider lax (or, alternatively,
pseudo) natural transformations $q_* F\to q_* G$. (For an explicit
description, substitute in the above diagrams by $1_{Ff}$ and $1_{Gf}$ the
idempotents $F1_f$ and $G1_f$, respectively, for any 1-cell $f$.) 
We call such a structure a {\em lax (or, alternatively, pseudo)
demitransformation} from $F$ to $G$.

This simplifies somewhat if $F$ and $G$ are in fact 2-functors:

\begin{proposition}
Let $F,G:\C\to\K$ be 2-functors. A lax demitransformation from $F$ to $G$
consists of a morphism $xC:FC\to GC$ in \K for each $C\in\C$ equipped 
with 2-cells $xf:Gf.xC\to xD.Ff$ satisfying conditions 
\eqref{eq:LN0} and \eqref{eq:LN1}
%(LN0) and (LN1)
only.  
\end{proposition}

\proof
We shall show that $\overline{x}C$ is just $x1_C$, and that all conditions 
involving it are then automatically satisfied. First of all, by compatibility 
with composition \eqref{eq:LN1}
%(LN1),
$x1_C$ is clearly idempotent. By \eqref{eq:DLN2} we have 
$x1_C.\overline{x}C=\overline{x}C$, while the fact that $x1_C$ is a 2-cell
in $Q_*\K$ gives $x1_C.\overline{x}C=x1_C$. Thus $\overline{x}C$ is 
necessarily just $x1_C$.

So now we define $\overline{x}C$ to be $x1_C$. Clearly 
\eqref{eq:LN0}, \eqref{eq:LN1}, and 
%(LN0), (LN1), and
\eqref{eq:DLN2} hold; we need only check that the composites
$$\xymatrix @R1pc {
Gf.xC \ar[r]^{xf} & xD.Ff \ar[r]^{x1_D.1_{Ff}} & xD.Ff}
\quad\textrm{and}\quad
\xymatrix @R1pc {
Gf.xC \ar[r]^{1_{Gf}.x1_C} & Gf.xC \ar[r]^{xf} & xD.Ff}$$
are both $xf$; these are both instances of compatibility with 
composition \eqref{eq:LN1}. %(LN1).
\endproof

\subsection{Modifications}
Finally we consider morphisms between lax natural transformations,
called {\em modifications}. 
In the case of lax demitransformations
$(x,\overline{x}),(y,\overline{y}):F\to G$, we retain the same word: a
modification from $(x,\overline{x})$ to $(y,\overline{y})$ consists 
of a 2-cell $\xi C:xC\to yC$ in \K for each $C\in\C$, subject to the 
usual condition expressed by commutativity of \eqref{eq:M} %(M), 
as well as
the extra condition represented by \eqref{eq:DM}: %(DM): 

\begin{minipage}{0.5\linewidth}
  \begin{equation}
    \label{eq:M}\tag{M}
\xymatrix{
Gf.xC \ar[r]^{xf} \ar[d]_{1_{Gf}.\xi C} & xD.Ff \ar[d]^{\xi D.1_{Ff}} \\
Gf.yC \ar[r]_{yf} & yD.Ff  }    
  \end{equation}
\end{minipage}
\begin{minipage}{0.5\linewidth}
  \begin{equation}
    \label{eq:DM}\tag{DM}
\xymatrix{
xC \ar[r]^{\xi C} \ar[dr]^{\xi C} \ar[d]_{\overline{x}C} &
yC \ar[d]^{\overline{y}C} \\
xC \ar[r]_{\xi C} & yC 
}  
  \end{equation}
\end{minipage}

% $$
% \textrm{(M)}
% \qquad\qquad
% \textrm{(DM)}
% $$
\noindent
Of course in the case of lax natural transformations, $\overline{x}C$
and $\overline{y}C$ are identities and so commutativity of \eqref{eq:DM} % (DM) 
is automatic.

\section{Limits}\label{sect:limits}

We now turn to the notion of limit within our ``weak world''. Because
of the well-established sense of ``weak limit'', referred to in the
introduction, we henceforth drop completely the epithet ``weak'', and speak
only of {\em demilimits}.

For an ordinary functor $S:\C\to\K$, the limit of $S$, if it exists,
is defined as the representing object of the functor from \C to \Set
sending $C\in\C$ to the set of cones under $S$ with vertex $C$.
Such a cone is of course just a natural transformation from the constant
functor $\Delta C$ at $C$ to $S$. Thus the notion of limit depends, among
other things, on the notion of naturality. In the 2-categorical context,
there is the possibility to replace naturality by lax naturality, giving
rise to a notion of lax limit \cite{Street-Cat-limits}; but in light of the
previous section we could instead consider (lax) deminaturality and so obtain
a notion of demilimit. (For more on bilimits and lax limits see
\cite{Street-Cat-limits} or \cite{Kelly-limits}.) 

For (small) 2-categories $\C$ and $\K$ we write $[\C,\K]$ for the usual
2-category of 2-functors from $\C$ to $\K$, with 2-natural transformations as
1-cells and modifications of 2-cells. Clearly, there is a natural bijection
between 2-functors $\K\to [\C,\K]$ and $\K\times \C\to \K$.
We write 
$\Ps(\C,\K)\lax$ for the bicategory of lax functors from $\C$ to $\K$,
with pseudonatural
transformations as 1-cells, and modifications as 2-cells. We denote by $J$
the fully faithful inclusion $\Catcc\to \Cat$, of the full sub-2-category of
\Cat consisting of the Cauchy complete categories. 

\subsection{Weighted bilimits}

Let $\C$ be a 2-category and $S:\C\to\K$ be a lax
demifunctor (of course this includes the case of an ordinary lax functor, or
indeed of a 2-functor); let $F:\C\to\Cat_{cc}$ be a 2-functor. The {\em
demilimit of $S$ weighted by $F$} is defined to be the $JF$-weighted bilimit
of the lax functor $S:\C\to Q_*\K$.  
That is, an object $\dl FS$ of $Q_*\K$ (i.e. of \K) equipped with a
pseudonatural equivalence  
$$
Q_*\K(-,\dl FS)\simeq \Ps(\C,\Cat)\lax(JF,JQ_*\K(-,S)),
$$
equivalently, 
$$
Q_*\K(-,\dl FS)\simeq
\Ps(\C,\Cat_{cc})\lax(F,Q_*\K(-,S)).
$$
If also the domain 2-category $\C$ is locally Cauchy complete, then this
notion of demilimit of a lax demifunctor $\C\to \K$ coincides with the bilimit
of the respective lax functor $\C\to Q_*\K$ in the $\Cat_{cc}$ enriched sense.

Similarly we have the {\em lax demilimit $\mathrm{dll}(F,S)$} defined by
$$Q_*\K(-,\mathrm{dll}(F,S))\simeq\Lax(\C,\Cat_{cc})\lax(F,Q_*\K(-,S)).$$
Note that
$\mathrm{dll}(F,S)$ can be constructed as $\dl{F'}S$ in terms of an
appropriate weight $F'$. 

A demi{\em colimit} in \K is of course just a demilimit notion in $\K\op$.

\begin{example}\label{ex:diag}
Let $T:\C\to \Catcc$ be the 2-functor constant on the terminal category
$1$. The demilimit $\dl TS$ of a lax demifunctor $S:\C\to \K$ is defined by
the pseudonatural equivalence 
$$
Q_*\K(-,\dl TS)\simeq \Ps(\C,\Cat)\lax(JT,JQ_*\K(-,S))
\simeq  \Ps(\C,Q_* \K)\lax (\Delta(-),S),
$$ 
cf. \cite[Section 4]{Street-Cat-limits}, where $\Delta:\K\to [\C,\K] \to
\mathrm{Ps}(\C,\K)\lax$ denotes the diagonal 2-functor (with the first arrow
corresponding to the first projection $\K\times \C\to \K$ and the second one
being the obvious inclusion). Thus for this particular weight $T$, the
demilimit $\dl TS$ is directly related to the bicategory of lax demifunctors,
demitransformations and modifications in Section \ref{sect:Ps}.
\end{example}

Our primary focus will be the case of Eilenberg-Moore objects, to which
we turn in the following section.

\subsection{Eilenberg-Moore objects}\label{sect:EM}

We have already seen the notion of Eilenberg-Moore object: for 
a monad $t$ on a category $B$ we write $B^t$ for the category of 
$t$-algebras; for a monad $t$ on an object $B$ of a 2-category \K, 
we write $B^t$ for the representing object in a representation
$$\K(A,B^t) \cong \K(A,B)^{\K(A,t)}$$
where $\K(A,t)$ is the induced monad on the hom-category $\K(A,B)$.
As a first generalization, we do not ask for an isomorphism, but just
a pseudonatural equivalence
$$\K(A,B^t) \simeq \K(A,B)^{\K(A,t)}$$
and we then call $B^t$ a {\em bicategorical EM-object}.

This fits into the framework of weighted limits of the previous section, more
specifically the situation in Example \ref{ex:diag}. 
We take \C to be the terminal 2-category, then a monad $(B,t)$ in \K is 
simply a lax functor $S:\C\to\K$. 
We take the weight $T:\C\to\Cat$ to be the  2-functor constant at the
terminal category. Then a lax natural transformation from $T$ to $\K(A,S)$
consists of a single component $b:A\to B$, with lax naturality constraint in
the form of a 2-cell  $\beta:tb\to b$, with the conditions stating $(b,\beta)$
is a $t$-algebra. 
Thus the bicategorical Eilenberg-Moore object of $t$ is just the
lax bilimit of the corresponding lax functor $S:1\to\K$ (weighted
by $T:1\to\Cat$). 

We now look in detail at the ``demi'' version. A demimonad in \K is
a monad in $Q_*\K$. For such a monad, a bicategorical EM-object
amounts to a representing object for $Q_*\K(A,B)^{Q_*\K(A,t)}$, as a 2-functor 
of $A\in Q_*\K$. Then we seek a pseudonatural equivalence 
\begin{equation}\label{eq:bicEM}
Q_*\K(A,B^{(t)})\simeq Q_*\K(A,B)^{Q_*\K(A,t)}.
\end{equation}
The right hand side is the category of demiactions of $t$ 
(cf. Proposition \ref{prop:demiaction}). 
The universal property guarantees a morphism $u:B^{(t)}\to B$ with 
a demiaction, $\psi:tu\to u$, such that for any demiaction $(a:A\to
B,\alpha:ta\to a)$, there exists a morphism $(x,\overline x):A\to B^{(t)}$ and
an isomorphism of demialgebras $\xi:(ux,\psi \overline x) \cong(a,\alpha)$.   

This becomes particularly simple in the case where $t$ is actually
a monad in \K. Then the right hand side becomes just $\K(A,B)^{(\K(A,t))}$,
and we seek a pseudonatural equivalence 
$$
Q_*\K(A,B^{(t)})\simeq \K(A,B)^{(\K(A,t))}.
$$
We can make this more explicit as follows. There is a morphism
$u:B^{(t)}\to B$, equipped with an action $\psi:tu\to u$
satisfying the associative law $\psi.t\psi=\psi.\mu
u$, but not required to satisfy the unit law. From the requirement that
it induces an equivalence, it has the following universal properties. 
By essential surjectivity on objects, for any morphism $a:A\to B$ and any
demiaction $\alpha:ta\to a$, there exists a morphism $(x,\overline x):A\to
B^{(t)}$ in $Q_*\K$ and an isomorphism $\xi:(u,
\psi_1:=\psi.\eta u)(x,\overline x)\cong (a,\alpha.\eta a)$
in $Q_*\K$ for which the diagram  
$$\xymatrix @R.7pc {
tux \ar[r]^{t\xi} \ar[ddd]_{\psi x} & ta \ar[ddd]^{\alpha} \\
\\
\\
ux \ar[r]_{\xi} & a 
}$$
commutes. Furthermore, there is a 2-dimensional aspect coming from the fact
that $(u,\psi)$ induces a fully faithful functor. 
Let $(x,\overline x),(y,\overline y):A\to B^{(t)}$ be given. For any
$\zeta:ux\to uy$ for which the diagrams
$$\xymatrix @R=1pc{
tux \ar[r]^{t\zeta} \ar[dd]_{\psi x} & tuy
\ar[dd]^{\psi y} && 
ux \ar[r]^{\zeta} \ar[d]_{\psi_1 x}
 \ar[ddr]^{\zeta} & 
uy \ar[d]^{\psi_1 y}
\\
&&& ux\ar[d]_-{u{\overline x}}&uy\ar[d]^-{u{\overline
    y}}\\ 
ux \ar[r]_{\zeta} & uy && ux \ar[r]_{\zeta} & uy }$$
commute, there is a unique 2-cell $\zeta':x\to y$ with
$\zeta'{\overline   x}=\zeta'= {\overline y}\zeta'$ and 
$u\zeta'.\psi_1 x=\zeta$. 

For any demimonad $(B,t)$ in a 2-category $\K$,
the image of the object $((t,\mu_1),\mu)$ of $Q_*\K(B,B)^{Q_*\K(B,t)}$ under
the isomorphism \eqref{eq:bicEM} provides a left adjoint
$(f,\overline{f}):B\to B^{(t)}$ of the 1-cell $(u,\psi_1)$
in $Q_*\K$ above. By virtue of the
universal property we have seen, the corresponding monad in $Q_*\K$ 
is isomorphic to $(B,t)$. This means the existence of 2-cells
$\chi:uf \to t$ and $\chi':t\to uf$
obeying the normalization conditions
$$
\xymatrix@C=20pt{
uf\ar[r]^-{\chi}\ar[d]_-{\psi_1\overline{f}} 
\ar[rd]^-{\!\!\chi}&
t \ar[d]^-{\mu_1}&&
t\ar[d]_-{\mu_1}\ar[rd]^-{\!\!\chi'}\ar[r]^-{\chi'}&
uf\ar[d]^-{\psi_1\overline{f}}&&
uf\ar[rd]_-{\psi_1\overline{f}}
\ar[r]^-{\chi} &
t\ar[d]^-{\chi'}&&
t\ar[rd]_-{\mu_1}\ar[r]^-{\chi'}&
uf\ar[d]^-{\chi}\\
uf\ar[r]_-{\chi}&t&&
t\ar[r]_-{\chi'}&
uf&&
&uf&&
&t
}
$$
and the `$t$-linearity' condition
$$
\xymatrix{
tuf\ar[r]^-{t\chi} \ar[d]_-{\psi f}&
t^2\ar[d]^-\mu\\
uf\ar[r]_-{\chi} &t}
\quad \textrm{equivalently,}\quad
\xymatrix{
t^2\ar[d]_-\mu\ar[r]^-{t\chi'} &
tuf\ar[d]^-{\psi f}\\
t\ar[r]_-{\chi'}&
uf.}
$$
The counit $(f,\overline{f})(u,\psi_1)\to 1$ takes the form of a map 
$\epsilon:fu\to 1$ for which the diagrams
$$\xymatrix{
ufu \ar[r]^{\chi u} \ar[d]_{\psi_1 fu} & tu \ar[d]^{\psi} &
fu \ar[r]^{\overline{f}\psi_1} \ar[dr]_{\epsilon} & fu \ar[d]^{\epsilon} \\
ufu \ar[r]_{u\epsilon} & u && 1 }$$
commute.

As usual \cite{ftm} there are various dualities. We write 
$\K\co$ for the 2-category obtained from \K by formally reversing the
direction of the 2-cells, but leaving the 1-cells unchanged. A monad
in $\K\co$ is a comonad in \K, and its demi-EM-object is just called
the demi-EM-object of the comonad. We write $\K\op$ for the 2-category
obtained from \K by formally reversing the direction of the 1-cells, but
leaving the 2-cells unchanged. A monad in $\K\op$ is still just a monad,
but the demi-EM-object is now a {\em colimit} in \K, called the 
demi-KL-object (KL for Kleisli).  
Finally we write $\K\coop$ for the 2-category obtained from \K by 
reversing both the 1-cells and the 2-cells. A monad in $\K\coop$ is a
comonad in \K; its demi-EM-object is called the demi-KL-object of the comonad.

\section{Free completions}

For a small category \C, the presheaf category $[\C\op,\Set]$ is the
free completion of \C under colimits. More generally, the free completion
of \C under some class of colimits is the closure of the representables
in $[\C\op,\Set]$ under those colimits. For example, the free completion 
of \C under coproducts is the full subcategory of $[\C\op,\Set]$ consisting
of those objects which are coproducts of representables. 

Furthermore, this remains true
in the enriched context: if \V is a complete and cocomplete symmetric monoidal
closed category, and \C is a small \V-category, the presheaf  
category $[\C\op,\V]$ is the free completion of \C under colimits; and 
the closure of the representables in $[\C\op,\V]$ under a given class
of colimits is the free completion under those colimits \cite{Kelly-book}. 
In particular this can be done in the case $\V=\Cat$ of 2-categories, leading
to a description of the free completion of a 2-category under Kleisli objects,
or dually under Eilenberg-Moore objects: this was the basis for the main 
construction in \cite{ftm2}. A potentially tricky aspect of these free
completions, is that one does not know how many steps may be involved in
forming the closure of the representables under some class of colimits: after
each step there will be new diagrams of which to form the colimit, and this
process could potentially continue transfinitely. In the case of Kleisli
objects, however, it terminates after a single step; this is basically because
a functor $f:C\to D$ exhibits $D$ as a Kleisli object if and only if it is
bijective on objects and has a right adjoint, and such functors are closed
under composition. 

We shall now consider ``demi'' versions of these ideas. In fact we treat
in detail only the case of completions under demi-KL-objects, but many 
other classes of demicolimits can be handled in similar fashion. 
In Section \ref{sect:limits} we defined demi-KL-objects as bilimits with
respect to a certain $\Catcc$-valued weight. When taking the free completion
under these bilimits, the key idea is to work not with categories enriched in
\Cat (2-categories) but rather with categories enriched in \Catcc, the full
2-subcategory of \Cat consisting of the Cauchy complete categories.  

The category \Catcc is Cartesian closed, so there is no problem enriching
over it:  a \Catcc-category is precisely a \Cat-category in which idempotent
2-cells split. The problem is that, as a category, \Catcc is neither complete
nor cocomplete, and so we cannot apply the Kelly theorem. One way around 
this is to note that although \Catcc is not complete or cocomplete as a
category, it is bicategorically complete and cocomplete (as a 2-category).
We can therefore use a bicategorical variant of the Kelly theorem, which we
prove in the Appendix. 
For a small \Catcc-category \C, we write $\Hom(\C\op,\Catcc)$ for the
2-category (in fact \Catcc-category) of pseudofunctors, pseudonatural
transformations, and modifications from $\C\op$ to \Catcc. Then
$\Hom(\C\op,\Catcc)$ is the free completion of the \Catcc-category \C under
bicategorical \Catcc-colimits, while the free completion under bicategorical
KL-objects is the closure under such of the representables in
$\Hom(\C\op,\Catcc)$. 

Once again, this process is potentially transfinite, but 
just as in the case of completion under ordinary
Kleisli objects, considered in \cite{ftm2}, the process terminates after 
a single step. 
The key observation here is the following lemma. Before stating it, it is
useful to define a functor $f:A\to B$ to be {\em quasi-surjective on objects}
if every object $b\in B$ is a retract of some $fa$ with $a\in A$. Such
functors are clearly closed under composition.

\begin{lemma}
A morphism $f:A\to B$ is of bicategorical Kleisli type in \Catcc if and 
only if it has a right adjoint and is quasi-surjective on objects.
\end{lemma}

\proof
Let $A$ be a Cauchy complete category, $t$ a monad on $A$, and $f_t:A\to A_t$
its Kleisli category; this is also a bicategorical Kleisli object in \Cat. 
Since $Q:\Cat\to\Catcc$ is left biadjoint to the inclusion $\Catcc\to\Cat$,
it preserves bicategorical Kleisli objects, and so the bicategorical Kleisli
object in \Catcc of $t$ is the composite
$$\xymatrix{
A \ar[r]^{f_t} & A_t \ar[r]^-{q} & QA_t. }$$
Now $A$ is Cauchy complete, and limits in the Eilenberg-Moore category $A^t$
can be formed as in $A$, so $A^t$ is also Cauchy complete. The canonical 
comparison $A_t\to A^t$ is fully faithful, and so $QA_t$ can be constructed,
up to equivalence, as the full subcategory of $A^t$ consisting of all
retracts of free algebras. It follows that the composite
$$\xymatrix{
QA_t \ar[r]^{j} & A^t \ar[r]^{u^t} & A }$$
where $j$ is the inclusion, is right adjoint to $qf_t$. On the other hand
$qf_t$ is clearly quasi-surjective on objects. This proves one half of 
the characterization.

Suppose conversely that a functor $f:A\to B$ in \Catcc has a right adjoint
$f\dashv u$ and is quasi-surjective on objects. We may form the induced 
monad $t$ on $A$; then the Kleisli category $A_t$ can be constructed by 
factorizing $f$ as an identity on object functor $f:A\to A_t$ followed
by a fully faithful one $q:A_t\to B$. Now $B$ is Cauchy complete,
and contains $A_t$ as a full subcategory, while every object of $B$ is
a retract of one in $A_t$; it follows that $B$ is equivalent to the Cauchy
completion $QA_t$ of $A_t$.
\endproof

Since bicategorical colimits in $\Hom(\C\op,\Catcc)$ are constructed pointwise,
we get a corresponding characterization of morphisms in $\Hom(\C\op,\Catcc)$
which are of bicategorical Kleisli type: the pseudonatural transformations
which pointwise have right adjoints and are quasisurjective on objects. Once
again, such morphisms are clearly closed under composition. It is this last
fact which means that we need only consider bicategorical Kleisli objects of  
monads on representables. 

At this point we may simply write down an explicit description of the 
free $\Catcc$-completion $\KL\dm(\K)$ of a small \Catcc-category
\C under bicategorical Kleisli objects. An object is a monad $(A,t)$ in \C
(with multiplication $\mu$  and unit $\eta$ understood). This 
generates a monad in $\Hom(\C\op,\Catcc)$ on the representable $\C(-,A)$. 
The Kleisli object is formed by first constructing the pointwise Kleisli
object in \Cat, then applying $Q$, to get 
$$\xymatrix{
\C(X,A) \ar[r]^-{F} & \C(X,A)_{\C(X,t)} \ar[r]^-{q} & Q(\C(X,A)_{\C(X,t)}). }$$
A morphism from $(A,t)$ to $(B,s)$ should be a pseudonatural transformation
$$\xymatrix{
Q(\C(X,A)_{\C(X,t)}) \ar[r] & Q(\C(X,B)_{\C(X,s)}) }$$
with values in \Catcc, or equivalently a pseudonatural transformation
$$\xymatrix{
\C(X,A)_{\C(X,t)} \ar[r] & Q(\C(X,B)_{\C(X,s)}) }$$
with values in \Cat, which in turn amounts to a pseudonatural transformation
$$\xymatrix{
\C(X,A) \ar[r] & Q(\C(X,B)_{\C(X,s)}) }$$
equipped with an op-action of $\C(X,t)$. By Yoneda the pseudonatural
transformation amounts to an object of $Q(\C(A,B)_{\C(A,s)})$; that is,
a morphism $f:A\to B$ equipped with a 2-cell $\phi_1:f\to sf$ which
is idempotent in $\C(A,B)_{\C(A,s)}$, or equivalently which satisfies 
$\mu f.s\phi_1.\phi_1=\phi_1$.
The op-action consists of a morphism in $Q(\C(A,B)_{\C(A,s)})$ from 
$(ft,\phi_1 t)$ to $(f,\phi_1)$, satisfying associativity and unitality
conditions. This then amounts to a 2-cell $\phi:ft\to sf$ in \K 
satisfying in addition to associativity and unitality two further
normalization conditions. The unitality condition says that $\phi.f\eta$ is
just $\phi_1$; it then turns out that the normalization conditions follow from
the single associativity condition $\mu f.s\phi.\phi
t=\phi.f\mu$. (Idempotency of $\phi_1$ is then automatic.) 

To summarize the situation so far, an object $\KL\dm(\K)$  is a monad, such as
$(A,t)$. A morphism from $(A,t)$ to $(B,s)$ is a morphism $f:A\to B$
in \K equipped with a 2-cell $\phi:ft\to sf$ satisfying the associativity
condition given above. What finally is a 2-cell between two 
such morphisms $(f,\phi)$ and $(g,\psi)$?
These should be modifications between the corresponding pseudonatural
transformations
$$\xymatrix{
Q(\C(X,A)_{\C(X,t)}) \ar[r] & Q(\C(X,B)_{\C(X,s)}) }$$
which reduce to modifications, compatible with the op-actions of $t$,  between
pseudonatural transformations 
$$\xymatrix{
\C(X,A) \ar[r] & Q(\C(X,B)_{\C(X,s)}) }$$
which by Yoneda amount to 2-cells $\rho:f\to sg$ subject to two conditions
stated in the theorem below; one gives compatibility with the op-actions, the
other is a normalization condition.

\begin{theorem}
Let \K be  a 2-category in which idempotent 2-cells split. The
free completion of \K as \Catcc-category under bicategorical Kleisli
objects, or equivalently the free completion of
\K under demi-KL-objects, is the evident 2-category $\KL\dm(\K)$ in which
\begin{enumerate}[(i)]
\item an object is a monad $(A,t)$ in \K;
\item a morphism from $(A,t)$ to $(B,s)$ is a 1-cell $f:A\to B$ in \K
  equipped with a 2-cell $\phi:ft\to sf$ for which the following diagram
  commutes; 
$$\xymatrix{
ftt \ar[r]^{\phi t} \ar[d]_{f\mu} & sft \ar[r]^{s\phi} & ssf \ar[d]^{\mu f} \\
ft \ar[rr]_{\phi} && sf }$$
\item a 2-cell from $(f,\phi)$ to $(g,\psi)$ is a 2-cell $\rho:f\to sg$
for which the following diagrams commute.
$$\xymatrix{
ft \ar[r]^{\rho t} \ar[d]_{\phi} & sgt \ar[r]^{s\psi} & ssg \ar[d]^{\mu g} &&
   sg \ar[r]^{sg\eta} & sgt \ar[r]^{s\psi} & ssg \ar[d]^{\mu g} \\
sf \ar[r]_{s\rho} & ssg \ar[r]_{\mu g} & sg && f \ar[u]^{\rho} \ar[rr]_{\rho} &&
sg }$$
\end{enumerate}
\end{theorem}

There is a formal dual of this, involving EM- rather than KL-objects. 
We write $\EM\dm(\K)$ for $\KL\dm(\K\op)\op$. This is exactly the
2-category $\EM^w(\K)$ of \cite{Bohm:wftm}. 

\begin{corollary}\label{cor:wftm}
If \K is a 2-category in which idempotent 2-cells split, then $\EM\dm(\K)$
is the free $\Catcc$-completion of \K under demi-EM-objects. 
\end{corollary}

What about the case of a general 2-category \K? There is a forgetful 2-functor 
from the 2-category of \Catcc-categories with demi-KL-objects to the 2-category
of 2-categories, and this forgetful 2-functor has a left biadjoint whose
object map can be constructed by first applying $Q_*$, then the construction
given above. We write $\KL\dm(\K)$ for the \Catcc-category obtained by
applying this left biadjoint to a 2-category \K, and call it the free
\Catcc-category with demi-KL-objects on \K. An object of $\KL\dm(\K)$ is just
a demimonad in \K; we write this as $(A,t)$, with remaining structure
$(\mu_2,\mu_1,\mu_0)$ omitted from the notation. A 1-cell from $(A,t)$ to
$(B,s)$ consists of a 1-cell $(f,\bar{f})$ in $Q_*\K$ equipped with a 2-cell
$\phi:(f,\bar{f})(t,\mu_1)\to(s,\mu_1)(f,\bar{f})$ in $Q_*\K$ satisfying
associativity. (We shall see shortly that a simplification is possible.) 
A 2-cell from $(f,\bar{f},\phi)$ to $(g,\bar{g},\psi)$ is a 2-cell
$\rho:(f,\bar{f})\to(s,\mu_1)(g,\bar{g})$ satisfying the two conditions
above. A 2-cell $(f,\bar{f})\to(s,\mu_1)(g,\bar{g})$ is a 2-cell 
$\rho:f\to sg$ such that $\rho\bar{f}=\rho=\mu_1g.s\bar{g}.\rho$; the 
other two conditions are unchanged.

Consider a 1-cell $(f,\bar{f},\phi):(A,t)\to(B,s)$. 
Let $\phi_1=\phi.f\eta:f\to sf$. This clearly defines a 
2-cell from $(f,\bar{f})\to(s,\mu_1)(f,1)$,
and compatibility with the op-action holds by 
$\mu f.s\phi_1.\phi=\mu f.s\phi.sf\eta.\phi=\mu f.s\phi.\phi
t.ft\eta=\phi.f\mu.ft\eta=\phi.f\mu_1=\phi.f\mu.f\eta t=\mu f.s\phi.\phi
t.f\eta t=\mu f.s\phi.\phi_1 t$ and finally the normalization condition 
by idempotency of $\varphi_1$; i.e.
$\mu f.s\phi.sf\eta.\phi_1=\mu f.s\phi.sf\eta.\phi.f\eta=\mu f.s\phi.\phi
t.ft\eta.f\eta=\phi.f\mu.ft\eta.f\eta=\phi.f\eta=\phi_1$; thus $\phi_1$ 
is a 2-cell from $(f,\bar{f},\phi)$ to $(f,1,\phi)$.

Similarly, $\phi_1$ is clearly a 2-cell from $(f,1)\to(s,\mu_1)(f,\bar{f})$,
and compatibility with the op-actions and the normalization condition hold
exactly as before, so we have a 2-cell from $(f,1,\phi)$ to
$(f,\bar{f},\phi)$, clearly inverse to the previous one.  

Thus in our 1-cells, we may as well restrict to those of the form
$(f,1,\phi)$, which we henceforth write simply as $(f,\phi)$. This 
gives the following description of $\KL\dm(\K)$ for general \K:

\begin{theorem}
The free $\Catcc$-category with demi-KL-objects on a 2-category \K is 
the evident 2-category $\KL\dm(\K)$ in which
\begin{enumerate}[(i)]
\item an object is a demimonad $(A,t)$ in \K;
\item a morphism from $(A,t)$ to $(B,s)$ is a 1-cell $f:A\to B$ in \K
  equipped with a 2-cell $\phi:ft\to sf$ for which the following diagrams
  commute;
$$\xymatrix{
ftt \ar[r]^{\phi t} \ar[d]_{f\mu} & sft \ar[r]^{s\phi} & ssf \ar[d]^{\mu f} \\  
ft \ar[rr]_{\phi} && sf}
\hspace{2.87cm}
\xymatrix{
ft \ar[r]^{\phi} \ar[dr]^{\phi} \ar[d]_{f\mu_1} & sf \ar[d]^{\mu_1 f} \\
  ft \ar[r]_{\phi} & sf }$$
\item a 2-cell from $(f,\phi)$ to $(g,\psi)$ is a 2-cell $\rho:f\to sg$
for which the following diagrams commute.
$$\xymatrix{
ft \ar[r]^{\rho t} \ar[d]_{\phi} & sgt \ar[r]^{s\psi} & ssg \ar[d]^{\mu g} &&
   sg \ar[r]^{sg\eta} & sgt \ar[r]^{s\psi} & ssg \ar[d]^{\mu g} \\
sf \ar[r]_{s\rho} & ssg \ar[r]_{\mu g} & sg && f \ar[u]^{\rho} \ar[rr]_{\rho} &&
sg }$$
\end{enumerate}
\end{theorem}

The new condition on morphisms says that the following composites
$$\xymatrix{
ft \ar[r]^{\phi} & sf \ar[r]^{\eta sf} & s^2f \ar[r]^{\mu f} & sf}
\qquad \textrm{and}\qquad 
\xymatrix{
ft \ar[r]^{ft\eta} & ft^2 \ar[r]^{f\mu} & ft \ar[r]^{\phi} & sf
}$$
are both simply equal to $\phi$. This is not automatic, as can be seen
by taking $(A,t)=(B,s)$ and $f=1$: then the identity 2-cell on $1t\to t1$
does not satisfy this condition unless $(A,t)$ is actually a monad.

Once again, there is a dual result for demi-EM-objects:

\begin{theorem}\label{thm:EMdmK}
The free $\Catcc$-category with demi-EM-objects on a 2-category \K is the
evident 2-category $\EM\dm(\K)$ in which 
\begin{enumerate}[(i)]
\item an object is a demimonad $(A,t)$ in \K;
\item a morphism from $(A,t)$ to $(B,s)$ is a 1-cell $f:A\to B$ in \K
  equipped with a 2-cell $\phi:sf\to ft$ for which the following diagrams
  commute; 
$$\xymatrix{
ssf \ar[r]^{s\phi} \ar[d]_{\mu f} & sft \ar[r]^{\phi t} & ftt \ar[d]^{f\mu} \\
sf \ar[rr]_{\phi} && ft}
\hspace{2.85cm}
\xymatrix{
sf \ar[r]^{\phi} \ar[dr]^{\phi} \ar[d]_{\mu_1 f} & ft \ar[d]^{f\mu_1} \\
  sf \ar[r]_{\phi} & ft }$$
\item a 2-cell from $(f,\phi)$ to $(g,\psi)$ is a 2-cell $\rho:f\to gt$
for which the following diagrams commute.
$$\xymatrix{
sf \ar[r]^{s\rho} \ar[d]_{\phi} & sgt \ar[r]^{\psi t} & gtt \ar[d]^{g\mu} &&
   gt \ar[r]^{\eta gt} & sgt \ar[r]^{\psi t} & gtt \ar[d]^{g\mu} \\
ft \ar[r]_{\rho t} & gtt \ar[r]_{g\mu} & gt && f \ar[u]^{\rho} \ar[rr]_{\rho} &&
gt }$$
\end{enumerate}
\end{theorem}

\section{Formal theory of monads}\label{sec:ftm}

The basic ingredients of the formal theory of monads, as presented
in \cite{ftm}, are as follows. For any 2-category \K, there is a 
2-category $\Mnd(\K)$ whose objects are monads in \K, and a fully 
faithful 2-functor $\Id:\K\to\Mnd(\K)$, sending an object of \K to the
identity monad on that object. This 2-functor has a right adjoint 
if and only if \K has Eilenberg-Moore objects; the right adjoint then
takes a monad to its Eilenberg-Moore object. Furthermore, there is 
a monad \Mnd on the category \twocat of 2-categories and 2-functors, 
and the endofunctor part of \Mnd sends an object \K to $\Mnd(\K)$,
while $\Id:\K\to\Mnd(\K)$ is the component at \K of the unit of the monad.
An object of $\Mnd(\Mnd(\K))$ --- that is, a monad in $\Mnd(\K)$ --- is the
same thing as a distributive law, and the multiplication
$\Comp:\Mnd(\Mnd(\K))\to\Mnd(\K)$ of the monad \Mnd sends a distributive law
to the induced composite monad.  

In the sequel \cite{ftm2} to \cite{ftm}, a variant $\EM(\K)$ of $\Mnd(\K)$ was
proposed, with the same objects and 1-cells as $\Mnd(\K)$, but with a more
general notion of 2-cell.  
Once again, this is the object-part of a monad on \twocat, and the 
unit $\Id:\K\to\EM(\K)$ has a right adjoint if and only if \K has 
Eilenberg-Moore objects; but this time there is a conceptual explanation:
$\EM(\K)$ is the free completion of \K under Eilenberg-Moore objects. 
From this universal property of $\EM(\K)$, it follows immediately that
$\Id:\K\to\EM(\K)$ will have a right adjoint if and only if \K has
Eilenberg-Moore objects; in particular, since $\EM(\K)$ has Eilenberg-Moore
objects, we obtain the multiplication $\Comp:\EM(\EM(\K))\to\EM(\K)$. 

\subsection{Wreaths}
An object of $\EM(\EM(\K))$ -- that is, a monad in $\EM(\K)$ -- is more
general than a monad in $\Mnd(\K)$, because of the more general 2-cells in 
$\EM(\K)$. Thus we obtain a more general notion of distributive law, called a
{\em wreath} in \cite{ftm2}. 

When it comes to the weak version, we have in place of $\EM(\K)$ our
weak version $\EM\dm(\K)$, the free completion of \K under demi-EM-objects.
Once again, we can draw various immediate conclusions from this universal
property of $\EM\dm(\K)$; many of these were given concrete proofs in
\cite{Bohm:wftm}, using the concrete description of $\EM\dm(\K)$. For 
instance, writing once again $\Id_\K:\K\to\EM\dm(\K)$ for the inclusion we
have:

\begin{theorem}
For any 2-category $\K$, the inclusion $\Id_{Q_*\K}:Q_*\K\to\EM\dm(Q_*\K)\cong
\EM\dm(\K)$ has a right biadjoint if and only if $\K$ has demi-EM objects.
In particular, the inclusion $\Id_\K:\K\to\EM\dm(\K)$ has a right biadjoint
whenever \K has bicategorical EM-objects and idempotent 2-cells split.
\end{theorem}

\proof
Note that for any object $X$ and any demimonad $(A,t)$ in $\K$, both
categories $Q_*\K(X,A)^{Q_*\K(X,t)}\cong\Mnd(Q_*\K)((X,1),(A,t))$ and
  $\EM\dm(\K)((X,1),(A,t))$ are isomorphic. Hence the claim follows
from the definition of the demi-EM object via the pseudonatural equivalence
$Q_*\K(X,A^{(t)})\simeq Q_*\K(X,A)^{Q_*\K(X,t)}$.
\endproof

In particular, the locally Cauchy complete 2-category $\EM\dm(\K)$ does have
demi-EM-objects, and so the inclusion $\EM\dm(\K)\to\EM\dm(\EM\dm(\K))$ does
have a right biadjoint, which sends demimonads $((A,t),(s,\lambda))$ in
$\EM\dm(\K)$ to demimonads $(A,st)$ in \K. 
We might call a demimonad in $\EM\dm(\K)$ a {\em
demiwreath} in \K. As an instance of the preceding theorem, every demiwreath
induces a composite demimonad. This is a (minor) generalization of one
direction of \cite[Theorem~2.3]{Bohm:wftm}, in that it deals from the outset
with demimonads rather than monads. 

The demi-EM object of the composite demimonad $(A,st)$ is defined via the
pseudonatural equivalence
$$
Q_*\K(X,A^{(st)})\simeq \EM\dm(\K)((X,1),(A,st))\simeq
\EM\dm(\EM\dm(\K)) (((X,1),1), ((A,t),(s,\lambda))).
$$
On the other hand, as said above, whenever demi-EM objects exist in $\K$,
$\Id_{Q_*\K}:Q_*\K\to  \EM\dm(\K)$ has a right biadjoint $J$ sending
an object $(A,t)$ of $\EM\dm(\K)$ (i.e. demimonad in $\K$)  to the demi-EM
object $A^{(t)}$. It induces a pseudofunctor $\Mnd(J):\Mnd(\EM\dm(\K))\to
\Mnd(Q_*\K)$, taking a demimonad $((A,t),(s,\lambda))$ in $\EM\dm(\K)$ to
the demimonad $(J(A,t),J(s,\lambda))=(A^{(t)},J(s,\lambda))$ in $\K$. The
demi-EM object of this latter monad is defined via the pseudonatural
equivalence 
\begin{eqnarray*}
Q_*\K(X,(A^{(t)})^{(J(s,\lambda))})&\simeq&
\Mnd(Q_*\K)((X,1),\Mnd(J)((A,t),(s,\lambda)))\\
&\simeq&\Mnd(\EM\dm(\K))(((X,1),1),((A,t),(s,\lambda)))\\
&\cong&\EM\dm(\EM\dm(\K))(((X,1),1),((A,t),(s,\lambda))).
\end{eqnarray*}
Thus we conclude that, whenever demi-EM objects exist in $\K$, $A^{(st)}$ and
$(A^{(t)})^{(J(s,\lambda))}$ are equivalent objects of $Q_*\K$. This extends
some observations in \cite[Proposition 3.7]{Bohm:wftm}. 

\subsection{Lifting}
Another key aspect of the formal theory of monads is that, for a 2-category \K 
with Eilenberg-Moore objects, monad morphisms
from $(f,\phi):(A,t)\to(B,s)$ are in bijection with morphisms $f:A\to B$
equipped with liftings
$$\xymatrix{
A^t \ar[r]^{\overline{f}} \ar[d] & B^s \ar[d] \\
A \ar[r]_{f} & B}$$
of $f$ to the Eilenberg-Moore objects. (More generally this is true provided
that the Eilenberg-Moore objects $A^t$ and $B^s$ exist.)  Similarly, for two
such morphisms $(f,\phi),(g,\psi):(A,t)\to (B,s)$, a 2-cell $\rho:f\to g$
gives a 2-cell $(f,\phi)\to(g,\psi)$ in $\Mnd(\K)$ if and only if it lifts to
a 2-cell $\overline{f}\to\overline{g}$ between the corresponding lifted
morphisms from $A^t$ to $B^s$. On the other hand a 2-cell in $\EM(\K)$ from
$(f,\phi)$ to $(g,\psi)$ is just an arbitrary 2-cell 
$\overline{f}\to\overline{g}$. There are analogues of this for $\EM\dm(\K)$. 

In the previous section we have seen that, whenever demi-EM objects exist in
$\K$ (equivalently, bicategorical EM objects exist in $Q_*\K$),
$\Id_{Q_*\K}:Q_*\K\to \EM\dm(\K)$ possesses a right biadjoint $J$ with object
map $(B,t)\mapsto B^{(t)}$. The counit of the biadjunction is given by the
1-cell $(u,\psi):(B^{(t)},1)\to (B,t)$ from Section
\ref{sect:EM}, for any demimonad $(B,t)$, and the iso 2-cell 
\begin{equation}\label{eq:1-lift}
\xymatrix@C=35pt{
(A^{(s)},1)\ar[r]^-{(J(g,\lambda),1)}\ar[d]_-{(u ,\psi )}
\ar@{}[rd]|-{{\Huge\Downarrow}\  \xi}& 
(B^{(t)},1)\ar[d]^-{(u,\psi)}\\
(A,s)\ar[r]_-{(g,\lambda)}&(B,t)}
\end{equation}
for any demimonad morphism (i.e. 1-cell in $\EM\dm(\K)$) $(g,\lambda):(A,s)\to
(B,t)$. Explicitly, such an iso 2-cell is given by 2-cells
$\xi:uJ(g,\lambda)\to gu $ and $\xi':gu \to
uJ(g,\lambda)$ in $\K$ such that the normalization conditions
$$
\xymatrix@C=15pt{
uJ(g,\lambda)\ar[r]^-{\xi}\ar[rd]_-\xi&gu \ar[d]^-{\lambda_1}&
gu \ar[r]^-{\xi'}\ar[rd]_-{\xi'}&
uJ(g,\lambda)\ar[d]^-{\psi_1J(g,\lambda)} &
gu \ar[r]^-{\xi'}\ar[rd]_-{\lambda_1}&
uJ(g,\lambda)\ar[d]^-\xi&
uJ(g,\lambda)\ar[r]^-{\xi}\ar[rd]_-{\psi_1J(g,\lambda)}&
gu \ar[d]^-{\xi'}\\
&gu &
&uJ(g,\lambda)&
&gu &
&uJ(g,\lambda)
}
$$
and the `$t$-linearity' condition
$$
\xymatrix@R=10pt{
tuJ(g,\lambda)\ar[r]^-{t\xi}\ar[dd]_-{\psi J(g,\lambda)}&
tgu \ar[d]^-{\lambda u }\\
&gsu \ar[d]^-{g\psi }\\
uJ(g,\lambda)\ar[r]_-{\xi}&gu }
$$
hold, where we introduced the idempotent 2-cells $\psi_1:=\psi.\eta u$
and $\lambda_1:=g\psi .\lambda u . \eta g u $, corresponding to the
$t$-demialgebras $(u,\psi)$ and $(gu ,g\psi .\lambda u )$, respectively. Note
that in particular $g\psi _1.\xi=\xi=\xi.\psi_1J(g,\lambda)$ and 
$\xi'.g\psi _1=\xi'=\psi_1J(g,\lambda).\xi'$.

In what follows we show that the map $(g,\lambda)\mapsto J(g,\lambda)$
provides the object map of an equivalence between the hom category of some
2-category of monads; and an appropriately defined category of liftings
for demi-EM objects $A^{(s)}\to B^{(t)}$.
\begin{lemma}\label{lem:1-lift}
Consider demimonads $(A,s)$ and $(B,t)$ in a 2-category $\K$ in which demi-EM
objects exist. If for some 1-cells $g:A\to B$ and $h:A^{(s)}\to B^{(t)}$ there
exist 2-cells $\zeta:uh\to gu $ and $\zeta':gu \to uh$
such that $\zeta'.\zeta=\psi_1h$ and the normalization conditions 
$g\psi _1.\zeta=\zeta=\zeta.\psi_1h$ and
$\zeta'.g\psi _1=\zeta'=\psi_1h.\zeta'$ hold, 
then there exists a demimonad morphism
$(g,\lambda):(A,s)\to (B,t)$ such that $J(g,\lambda)$ is isomorphic to $h$.
\end{lemma}

\proof
The requested 1-cell $(g,\lambda):(A,s)\to (B,t)$ is constructed by
introducing $\lambda$ as the composite
$$
\xymatrix{
tg\ar[r]^-{tg\eta}&
tgs\ar[r]^-{tg\chi' }&
tgu f \ar[r]^-{t\zeta'f }& 
tuhf \ar[r]^-{\psi hf }&
uhf \ar[r]^-{\zeta f }&
gu f \ar[r]^-{g\chi }&gs},
$$
where the notations from Section \ref{sect:EM} are used. 
The corresponding idempotent 2-cell $\lambda_1:gu \to gu $ comes out
as $\zeta.\zeta'$. For the induced $t$-demiaction $g\psi .\lambda
u =\zeta.\psi h.t\zeta':tgu \to gu $, both $\zeta$ and $\zeta'$
are morphisms of $t$-demialgebras. Hence together with the canonical 2-cells
$\xi:uJ(g,\lambda)\to gu $ and $\xi':gu  \to
uJ(g,\lambda)$, they induce $t$-demialgebra morphisms $\xi'.\zeta:
uh\to uJ(g,\lambda)$ and $\zeta'.\xi:uJ(g,\lambda)\to
uh$. These are subject to the normalization conditions in Section
\ref{sect:EM} hence by universality of $(u,\psi)$ give rise to
mutually inverse isomorphisms $\alpha:h\to J(g,\lambda)$ and
$\alpha':J(g,\lambda)\to h$.
\endproof

By naturality of $\xi$, for any
2-cell $\varrho:(g,\lambda)\to (g',\lambda')$ in $\EM\dm(\K)$, the 2-cell
$J(\varrho): J(g,\lambda)\to J(g',\lambda')$ renders commutative
\begin{equation}\label{eq:J_2}
\xymatrix @R=10pt{
uJ(g,\lambda)\ar[dd]_-{uJ(\varrho)}\ar[r]^-\xi&
gu \ar[d]^-{\varrho u }\\
&g'su \ar[d]^-{g'\psi }\\
uJ(g',\lambda')\ar[r]_-{\xi}&
g'u }
\quad \textrm{equivalently,}\quad
\xymatrix @R=10pt{
gu \ar[d]_-{\varrho u }\ar[r]^-{\xi'}&
uJ(g,\lambda)\ar[dd]^-{uJ(\varrho)}\\
g'su \ar[d]_-{g'\psi }&\\
g'u  \ar[r]_-{\xi'}&uJ(g',\lambda')}
\end{equation}
The above two equivalent forms of the same equality provide us with two
symmetrical choices how to define a lifting of a 2-cell in $\K$ for demi-EM
objects: we can require either one to take a particularly simple
form.

\begin{lemma}\label{lem:lift}
Let $\K$ be a 2-category in which demi-EM objects exist and $(g,\lambda)$ and
$(g',\lambda')$ be demimonad morphisms $(A,s)\to (B,t)$ in $\K$.

\noindent
(1) For a 2-cell $\omega:g\to g'$, the following are equivalent.
\begin{enumerate}[(i)]
\item the following diagram commutes;
$$
\xymatrix{
tg\ar[rrrr]^-{t\omega}\ar[d]_-{\lambda}&&&&
tg'\ar[d]^-{\lambda'}\\
gs\ar[r]_-{\omega s}&
g's\ar[r]_-{\eta g's}&
tg's\ar[r]_-{\lambda' s}&
g's^2\ar[r]_-{g'\mu}&
g's}
$$
\item $\lambda'_1.\omega u $ is a $t$-demialgebra morphism
  $(gu ,g\psi .\lambda u )\to  (g'u ,g'\psi .\lambda'
  u )$; 
\item $\lambda'.\eta g'.\omega:g\to g's$ is a 2-cell $(g,\lambda)\to
  (g',\lambda')$ in $\EM\dm(\K)$;
\item there is a 2-cell $\plift \omega:J(g,\lambda)\to J(g',\lambda')$ such
  that the following diagram commutes.
$$
\xymatrix{
gu \ar[r]^-{\omega u }\ar[d]_-{\xi'}&
g'u \ar[d]^-{\xi'}\\
uJ(g,\lambda)\ar[r]_-{u\plift \omega}&
uJ(g',\lambda')}
$$
\end{enumerate}
If these assertions hold then $\plift \omega =J(\lambda'.\eta g'.\omega)$.

\noindent
(2) For a 2-cell $\omega:g\to g'$, the following are equivalent.
\begin{enumerate}[(i)]
\item the following diagram commutes;
$$
\xymatrix@R=12pt{
tg\ar[r]^-{t\eta g}\ar[dd]_-{\lambda}&
t^2g\ar[r]^-{t\lambda}&
tgs\ar[r]^{t\omega s}&
tg's\ar[d]^-{\lambda's}\\
&&&g's^2\ar[d]^-{g'\mu}\\
gs\ar[rrr]_-{\omega s}&&&
g's}
$$
\item $\omega u .\lambda_1$ is a $t$-demialgebra morphism
  $(gu ,g\psi .\lambda u )\to  (g'u ,g'\psi .\lambda'
  u )$; 
\item $\omega s.\lambda.\eta g:g\to g's$ is a 2-cell $(g,\lambda)\to
  (g',\lambda')$ in $\EM\dm(\K)$;
\item there is a 2-cell $\ilift \omega:J(g,\lambda)\to J(g',\lambda')$ such
  that the following diagram commutes.
$$
\xymatrix{
uJ(g,\lambda)\ar[r]^-{u\ilift \omega}\ar[d]_-{\xi}&
uJ(g',\lambda')\ar[d]^-{\xi}\\
gu \ar[r]_-{\omega u }&
g'u }
$$
\end{enumerate}
If these assertions hold then $\ilift \omega =J(\omega s.\lambda.\eta g)$.
\end{lemma}

\proof
Consider first part (1).
The diagram in part (i) is equivalent to $t$-linearity of the 2-cell in part
(ii). If this holds then the normalization conditions on the 2-cell in part
(ii) are automatic. Thus (i)$\Leftrightarrow$(ii). 
Similarly, the diagram in part (i) is equivalent to the first diagram in
Theorem \ref{thm:EMdmK} (iii) for taking the 2-cell in part (iii) as
``$\varrho$''. If this holds then the second condition in Theorem
\ref{thm:EMdmK} (iii) is automatic. Thus (i)$\Leftrightarrow$(iii). 
If the assertion in part (iii) holds, then we can obtain $\plift \omega$ in
part (iv) by applying the pseudofunctor $J$ to the 2-cell in part (iii). The
diagram in part (iv) is then just the second diagram in \eqref{eq:J_2} for the
2-cell in part (iii). 
Finally, assume that assertion (iv) holds. Then $\lambda'_1.\omega
u =\xi.\xi'.\omega u =\xi. u \plift \omega .\xi'$ is evidently
a morphism of $t$-demialgebras hence also (ii) holds.

Part (2) is proven similarly, using the first diagram in \eqref{eq:J_2}
instead of the second one.
\endproof

\begin{corollary}
Let $\K$ be a 2-category in which demi-EM objects exist and $(A,s)$ and
$(B,t)$ be demimonads in $\K$.

\noindent
(1) The following categories are equivalent.

\begin{enumerate}[(i)]
\item The category whose objects are quadruples $(g:A\to B,h:A^{(s)}\to
  B^{(t)},\zeta:uh\to gu ,\zeta':gu \to uh)$
such that $\zeta'.\zeta=\psi_1h$ and the normalization conditions 
$g\psi _1.\zeta=\zeta=\zeta.\psi_1h$ and
$\zeta'.g\psi _1=\zeta'=\psi_1h.\zeta'$ hold. Morphisms
$(g,h,\zeta,\zeta')\to (\tilde g,\tilde h,\tilde\zeta,\tilde\zeta')$ are pairs
$(\omega:g\to \tilde g,\varphi:h\to \tilde h)$ such that
$u\phi.\zeta'=\tilde \zeta'.\omega u $.
\item The category whose objects are demimonad morphisms $(A,s)\to (B,t)$ and
  morphisms $(g,\lambda)\to (g',\lambda')$ are 2-cells $\omega: g\to g'$
  rendering commutative the diagram in Lemma \ref{lem:lift} (1)(i).
\end{enumerate}

\noindent
(2) The following categories are equivalent.

\begin{enumerate}[(i)]
\item The category whose objects are quadruples $(g:A\to B,h:A^{(s)}\to
  B^{(t)},\zeta:uh\to gu ,\zeta':gu \to uh)$
such that $\zeta'.\zeta=\psi_1h$ and the normalization conditions 
$g\psi _1.\zeta=\zeta=\zeta.\psi_1h$ and
$\zeta'.g\psi _1=\zeta'=\psi_1h.\zeta'$ hold. Morphisms
$(g,h,\zeta,\zeta')\to (\tilde g,\tilde h,\tilde\zeta,\tilde\zeta')$ are pairs
$(\omega:g\to \tilde g,\varphi:h\to \tilde h)$ such that
$\tilde \zeta.u\phi=\omega u .\zeta$.
\item The category whose objects are demimonad morphisms $(A,s)\to (B,t)$ and
  morphisms $(g,\lambda)\to (g',\lambda')$ are 2-cells $\omega: g\to g'$
  rendering commutative the diagram in Lemma \ref{lem:lift} (2)(i).
\end{enumerate}
\end{corollary}

\proof
By Lemma \ref{lem:lift}, there are fully faithful functors from
the categories in parts (ii) to the respective categories in part (i). They
are essentially surjective on the objects by Lemma \ref{lem:1-lift}:
for any object $(g:A\to B,h:A^{(s)}\to B^{(t)},\zeta:uh\to
gu ,\zeta':gu \to uh)$ in part (i), the isomorphism $h\to
J(g,\lambda)$ in  Lemma \ref{lem:1-lift} and the identity morphism $g\to g$
constitute an isomorphism in the category in question.
\endproof

The categories in the above corollary are hom categories in evident
2-categories. Both parts (ii) amount to extensions of $\Mnd(\K)$ in two
inequivalent ways.  

As we recalled earlier, a distributive law is in fact a monad in
$\Mnd(\K)$. Weak distributive laws in \cite{street:wdl} are not the same as
monads in either generalization of $\Mnd(\K)$ in the above corollary. However,
following the lines in \cite{BLS:wdl}, they can be described as compatible
pairs of monads in both of them.

\appendix

\section{Free $\mathbf{Cat_{cc}}$-completions}\label{appendix}

The classical theory of weighted colimits and colimit completions can
be found in \cite{Kelly-book}. It applies for categories enriched in a
complete and cocomplete symmetric monoidal closed category \V. Here we
adapt it to deal with the case of categories enriched over the
Cartesian closed category \Catcc, which is neither complete not
cocomplete. It is, however, complete and cocomplete as a bicategory,
and this will be the basis of our approach.

We write $J:\Catcc\to\Cat$ for the fully faithful inclusion; it has a
left biadjoint $Q$. 
For 2-categories \A and \B  we write $[\A,\B]$ for the usual 2-category of
2-functors from \A to \B, with 2-natural transformations as 1-cells and
modifications as 2-cells. We write $\Hom(\A,\B)$ for the 2-category of
pseudofunctors, pseudonatural transformations, and modifications, and
$\Ps(\A,\B)$ for the 2-category of 2-functors, pseudonatural transformations,
and modifications.  

\begin{remark}
Let $F:\A\to\Catcc$ be a pseudofunctor. Then $JF:\A\to\Cat$ is also 
a pseudofunctor. It is pseudonaturally equivalent to a 2-functor
$G:\A\to\Cat$; but any category equivalent to a Cauchy complete one is
itself Cauchy complete, and so $G$ lands in \Catcc, and can be written
as $JH$ for some 2-functor $H:\A\to\Catcc$, which is then
pseudonaturally equivalent to $F$. Thus the 2-categories
$\Ps(\A\op,\Catcc)$ and $\Hom(\A\op,\Catcc)$ are biequivalent.
\end{remark}

Recall that for pseudofunctors $S:\D\to\K$ and $F:\D\op\to\Cat$, the
bicolimit $F*S$ is defined by a pseudonatural equivalence
$$\K(F*S,A) \simeq \Hom(\D\op,\Cat)(F,\K(S,A)).$$
(In fact it does no harm to suppose that the map from left to right is
strictly natural in $A$, and so is induced by a pseudonatural
$F\to\K(S,F*S)$, but the inverse equivalence, going from right to
left, will still only be pseudonatural.)

If \K is a \Catcc-category, then we may choose to restrict to the case
where \D is a \Catcc-category and $F$ lands in \Catcc.

\begin{proposition}
For a small \Catcc-category \A, the 2-category $\Hom(\A\op,\Catcc)$ is 
in fact a \Catcc-category with all bicolimits.
\end{proposition}

\proof
The existence of bicolimits follows from the fact that the fully faithful
inclusion of $\Hom(\A\op,\Catcc)$ in $\Hom(\A\op,\Cat)$ has a left biadjoint.

The idempotent splittings can be computed pointwise. 
\endproof

Let $\Phi$ be a class of \Catcc-weights, and \A a small
\Catcc-category. Write $\Phi(\A)$ for the closure in
$\Hom(\A\op,\Catcc)$ of the representables under
$\Phi$-bicolimits. (This can be formed as the intersection of all full
subcategories containing the representables and closed under
$\Phi$-bicolimits; since $\Hom(\A\op,\Catcc)$ is such a subcategory,
and the intersection of any collection of such subcategories is one,
the intersection clearly has the desired properties. It can also be
formed via a transfinite induction.)

We shall write $W:\Phi(\A)\to\Hom(\A\op,\Catcc)$ for the inclusion, and 
$Y:\A\to\Phi(\A)$ for the restricted Yoneda embedding. We wish to
prove that $\Phi(\A)$ is the free completion of \A under $\Phi$-bicolimits,
in the sense that for any \Catcc-category \K with $\Phi$-bicolimits,
composition with $Y$ induces a biequivalence
$$\PhiCoc(\Phi(\A),\K)\simeq\Hom(\A,\K)$$
of \Catcc-categories. 

Our first result holds by definition of $\Phi(\A)$:

\begin{proposition}
$\Phi(\A)$ has $\Phi$-bicolimits, preserved by $W$. \endproof
\end{proposition}

From now on we shall fix a \Catcc-category \K with $\Phi$-colimits.

\begin{proposition}
For any pseudofunctor $F:\A\to\K$, the (pointwise) left Kan extension
$\Lan_Y F:\Phi(\A)\to\K$ exists.
\end{proposition}

\proof
The formula for the pointwise left Kan extension is 
$$(\Lan_Y F)X = X*F$$
so we are to show that the bicolimit on the right exists for all
$X\in\Phi(\A)$. 

Let \B be the full sub-\Catcc-category of $\Phi(\A)$ consisting of
those objects $X$ for which $X*F$ does exist. Certainly \B contains
the representables, since $\A(-,A)*F$ can be taken to be $FA$. On the
other hand, if $S:\D\to\B$ and $\phi:\D\op\to\Catcc$ is in $\Phi$,
then each $SD*F$ exists, so we may write $S*F$ for the pseudofunctor 
sending $D$ to $SD*F$, and now $\phi*(S*F)$ exists, since \K has
$\Phi$-weighted bicolimits. But $\phi*(S*F)\simeq(\phi*S)*F$, and so 
$\phi*S$ is also in \B. Thus \B contains the representables and is
closed under $\Phi$-weighted bicolimits, so must be all of $\Phi(\A)$. 
\endproof

\begin{proposition}
In the setting of the previous proposition, $\Lan_Y F$ is $\Phi$-cocontinuous.  
\end{proposition}

\proof
Let $S:\D\to\Phi(\A)$ and $\phi:\D\op\to\Catcc$, with
$\phi\in\Phi$. We must show that $(\Lan_Y F)(\phi*S)\simeq\phi*(\Lan_Y
F)S$. But 
$$(\Lan_Y F)(\phi*S) = (\phi*S)*F\simeq\phi*(S*F)=\phi*(\Lan_Y F)S$$
hence the result.
\endproof

It now follows that the 2-functor 
$$\PhiCoc(\Phi(\A),\K)\to\Hom(\A,\K)$$
given by restriction along $Y$ has a left biadjoint given by left Kan
extension
along $Y:\A\to\Phi(\A)$. The component at $F:\A\to\K$ of the unit is
the canonical map $F\to\Lan_Y(F)Y$, which is an equivalence since $Y$
is one. (More concretely,
$\Lan_Y(F)YA=\Lan_Y(F)\A(-,A)\simeq\A(-,A)*F\simeq FA$, and so 
$\Lan_Y(F)Y\simeq F$.)
Thus it remains only to show that the counit is also an
equivalence, which amounts to 

\begin{proposition}
If $G:\Phi(\A)\to\K$ is $\Phi$-cocontinuous, then the canonical map 
$\Lan_Y(GY)\to G$ is an equivalence.
\end{proposition}

\proof
Let \B be the full subcategory of $\Phi(\A)$ consisting of those
objects $X$, for which $\Lan_Y(GY)X\to GX$ is an equivalence; in other
words,
for which $X*GY\to GX$ is an equivalence. Then $\A(-,A)*GY\simeq
GYA=G\A(-,A)$,
and so \B contains the representables. Suppose that $S:\D\to\Phi(\A)$
lands in \B, and that $\phi:\D\op\to\Catcc$ is in $\Phi$. Then
\begin{align}
G(\phi*S) &\simeq \phi*GS \tag{$G$ is $\Phi$-cocontinuous} \\
               &\simeq \phi*(S*GJ) \tag{$S$ lands in \B} \\
               &\simeq (\phi*S)*GJ \tag{associativity of $*$} 
\end{align}
and so $\phi*S$ is also in \B, and thus \B is closed under
$\Phi$-colimits.
Thus \B is all of $\Phi(\A)$.\endproof


\begin{thebibliography}{Bibliography}{}

\bibitem{AlAletal:w_Gal}
Alonso \'Alvarez, J. N.; Fern\'andez Vilaboa, J. M.; Gonz\'alez
Rodr\'{\i}guez, R.; Rodr\'{\i}guez Raposo, A. B. 
{\em Weak $C$-cleft extensions and weak Galois extensions.}  
J. Algebra  299  (2006),  no. 1, 276-293. 
MR2225776 (2007b:16081) 

\bibitem{AlAletal:w_smash}
Alonso \'Alvarez, J. N.; Fern\'andez Vilaboa, J. M.; Gonz\'alez
Rodr\'{\i}guez, R. 
{\em Weak Hopf algebras and weak Yang-Baxter operators.} 
J. Algebra 320 (2008), no. 5, 2101--2143. 
MR2437645 (2009g:16057) 

\bibitem{bicategories}
B\'enabou, Jean. 
{\em Introduction to bicategories.}  
[in:] Reports of the Midwest Category Seminar  pp. 1-77 Springer, Berlin, 1967.
MR0220789 (36 \#3841)

\bibitem{Bohm:wftm}
B\"ohm, Gabriella.
{\em The weak theory of monads.}  
Adv. Math. 225 (2010), no. 1, 1-32.
MR2669347

\bibitem{BLS:weak_bimonad}
B\"ohm, Gabriella; Lack, Stephen; Street, Ross.
{\em Weak bimonads and weak Hopf monads.}
J. Algebra 328 (2011), no. 1, 1-30.
MR2745551

\bibitem{BLS:wdl}
B\"ohm, Gabriella; Lack, Stephen; Street, Ross.
{\em On the 2-categories of weak distributive laws.}
Comm. Algebra (special volume dedicated to Mia Cohen), in press. 
Preprint available at {\tt http://arxiv.org/abs/1009.3454}.

\bibitem{BNSz:WHAI}
B\"ohm, Gabriella; Nill, Florian; Szlach\'anyi, Korn\'el. 
{\em Weak Hopf algebras. I. Integral theory and $C^*$-structure.}  
J. Algebra  221  (1999),  no. 2, 385-438. 
MR1726707 (2001a:16059) 

\bibitem{Brz:coring}
Brzezi\'nski, Tomasz. 
{\em The structure of corings: induction functors, Maschke-type
theorem, and Frobenius and Galois-type properties.}  
Algebr. Represent. Theory 5  (2002),  no. 4, 389-410.
MR1930970 (2003h:16060)

\bibitem{CaDeG:w_entw}
Caenepeel, Stefaan; De Groot, Erwin. 
{\em Modules over weak entwining structures.}  
[in:] New trends in Hopf algebra theory (La Falda, 1999), 
eds.: N. Andruskiewitsch, W.R. Ferrer Santos and H-J. Schneider. 
pp. 31-54, Contemp. Math., 267, Amer. Math. Soc., Providence, RI, 2000. 
MR1800705 (2001k:16069) 

%\bibitem{GoPoSt}
%Gordon, R.; Power, A. J.; Street, Ross.
%{\em Coherence for tricategories.}
%Mem. Amer. Math. Soc. 117 (1995), no. 558, vi+81 pp.
%MR1261589 (96j:18002)

\bibitem{Ha:face}
Hayashi, Takahiro. 
{\em Quantum group symmetry of partition functions of IRF models and its
application to Jones' index theory.}
Comm. Math. Phys.  157  (1993),  no. 2, 331-345.
MR1244870 (94h:82020)

\bibitem{Kelly-limits}
Kelly, Gregory Maxwell. 
{\em Elementary observations on $2$-categorical limits.}
Bull. Austral. Math. Soc.  39  (1989),  no. 2, 301-317.
MR0998024 (90f:18004) 

\bibitem{Kelly-book}
Kelly, Gregory Maxwell. 
{\em Basic concepts of enriched category theory.} 
London Mathematical Society Lecture Note Series, 64. Cambridge University
Press, Cambridge-New York, 1982. 245 pp. ISBN: 0-521-28702-2.
[Re-published in:] 
Reprints in Theory and Applications of Categories, 
No. 10 (2005) pp. 1-136.
MR0651714 (84e:18001)

\bibitem{www}
Lack, Stephen.
{\em Welcome to weak world.}
Talk at CT2009 conference in Cape Town, 2 July 2009.

\bibitem{ftm2}
Lack, Stephen; Street, Ross. 
{\em The formal theory of monads. II.} 
[in:] Special volume celebrating the 70th birthday of Professor Max Kelly.
J. Pure Appl. Algebra 175  (2002),  no. 1-3, 243-265.
MR1935981 (2003m:18007)

\bibitem{Lawvere-metric}
Lawvere, F. William.
{\em Metric spaces, generalized logic, and closed categories.}
Rend. Sem. Mat. Fis. Milano 43 (1973), 135-166.
MR0352214 (50 \#4701)

\bibitem{CWM}
Mac Lane, Saunders. 
{\em Categories for the working mathematician.} 
Second edition. Graduate Texts in Mathematics, 5. Springer-Verlag, New York,
1998. xii+314 pp. ISBN: 0-387-98403-8 
MR1712872 (2001j:18001)

\bibitem{Oc:q_grd}
Ocneanu, Adrian. 
{\em Quantized groups, string algebras and Galois theory for algebras.}  
[in:] Operator algebras and applications, 
eds.: D. E. Evans and M. Takesaki. 
Vol. 2,  pp. 119-172, London Math. Soc. Lecture Note Ser., 136, 
Cambridge Univ. Press, Cambridge, 1988 
MR0996454 (91k:46068)

\bibitem{PaSt:w_bimonad}
Pastro, Craig; Street, Ross. 
{\em Weak Hopf monoids in braided monoidal categories.}
Algebra Number Theory  3  (2009),  no. 2, 149-207. 
MR2491942 (2010i:18007)

\bibitem{ftm}
Street, Ross. 
{\em The formal theory of monads.}  
J. Pure Appl. Algebra  2  (1972), no. 2, 149-168. 
MR0299653 (45 \#8701)

\bibitem{Street-Cat-limits}
Street, Ross. 
{\em Limits indexed by category-valued $2$-functors.}  
J. Pure Appl. Algebra  8  (1976), no. 2, 149-181.
MR0401868 (53 \#5695) 

\bibitem{street:wdl}
Street, Ross. 
{\em Weak distributive laws.}
Theory Appl. Categ. 22 (2009) no. 12, 313-320.
MR2520974 (2010k:18004)

\bibitem{Szlach}
Szlach\'anyi, Korn\'el. 
{\em Adjointable monoidal functors and quantum groupoids.} 
[in:] “Hopf algebras in noncommutative geometry and physics”, Caenepeel, S.;
Van Oysaeyen, F. (eds.), pp. 291–307, Lecture Notes in Pure and Appl. Math.,
239, Dekker, New York, 2005.  
MR2106937 (2005m:16052)

\bibitem{Ya:g_Kac}
Yamanouchi, Takehiko. 
{\em Duality for generalized Kac algebras and a characterization of finite
groupoid algebras.}
J. Algebra  163  (1994),  no. 1, 9-50.
MR1257303 (95c:22010) 

\bibitem{Wis:q-adj}
Wisbauer, Robert.
{\em Adjunction contexts and regular quasi-monads.}
Preprint available at 
{\tt http://arxiv.org/abs/1101.1195}.
\end{thebibliography}
\end{document}